\documentclass[12pt]{amsart}

\usepackage{amsmath,amssymb}
\usepackage[all]{xy}

\newtheorem{theorem}{Theorem}[section]

\newtheorem{cor}[theorem]{Corollary}
\newtheorem{defn}[theorem]{Definition}
\newtheorem{prop}[theorem]{Proposition}
\newtheorem{pro}[theorem]{Program}
\newtheorem{lem}[theorem]{Lemma}
\newtheorem{exa}[theorem]{Example}

\newtheorem{rem}[theorem]{Remark}
\numberwithin{equation}{section}

\newcommand{\rd}[1]{\left\lfloor#1\right\rfloor}

\newcommand{\rar}{\rightarrow}

\newcommand{\ZZ}{\mathbb Z}
\newcommand{\CC}{\mathbb C}

\newcommand{\QQ}{\mathbb Q}

\newcommand{\PP}{\mathbb P}

\newcommand{\ep}{\varepsilon}

\newcommand{\sT}{\mathcal T}
\newcommand{\sB}{\mathcal B}

\newcommand{\sC}{\mathcal C}

\newcommand{\sV}{\mathcal V}
\newcommand{\sN}{\mathcal N}
\newcommand{\sF}{\mathcal F}
\newcommand{\sX}{\mathcal X}
\newcommand{\sY}{\mathcal Y}

\newcommand{\Oh}{\mathcal O}

\DeclareMathOperator{\wt}{wt}

\DeclareMathOperator{\Proj}{Proj}

\DeclareMathOperator{\Spec}{Spec}

\DeclareMathOperator{\ch}{Ch}
\DeclareMathOperator{\td}{Td}
\DeclareMathOperator{\Stac}{Stac}

\DeclareMathOperator{\Pro}{Proj}
\begin{document}
\date{}
\title{Explicit orbifold Riemann--Roch for quasismooth varieties} 

\author{Shengtian Zhou}
\begin{abstract} Considering quasismooth varieities as global $\CC^*$ quotients, we present a Riemann-Roch formula via general Riemann-Roch formula for quotient stacks. Furthermore, we give a parcing formula for Hilbert series associated to a polarized quasismooth projectively Gorenstein algebraic varieties with orbifold curves and dissident points, which is an extension of the result in \cite{BRZ}.

\end{abstract}

\maketitle
\tableofcontents
\section{Introduction}
A quasismooth variety is a variety with weighted projective space as ambient space with all its singularities inherited from its ambient space. Therefore  it has only the cyclic quotient singularities (orbifold loci) of possibly different dimensions. In \cite{BRZ} we presented a formula for parsing Hilbert series associated to a polarized quasismooth projectively Gorenstein variety with only isolated singularities and some attempt to do the same for $3$-folds. Here we want to generate this to higher dimensional varieties with possible orbifold loci of dimension $1$. We see that this type of parsing for Hilbert series are based on the Riemmann-Roch formula for respective varieties (\cite{R}, \cite{BS}). Even though there already exists Riemann-Roch theorem for singular varities, see for example \cite{BFM} for general treatment for singular varities and \cite{K} for V-manifolds, they are too abstract for our purpose. 

In this paper, we will first present in section \ref{Riemann-Roch} an explicit Riemann-Roch formula for quasismooth varieties. Following Nironi's treatment for weighted projective space in \cite{N}, we obtain this formula via the general formula for stacks of To\"{e}n \cite{To}. We will see that viewing quasismooth varieties as a quotient stack exposes their structures in a transparent way.  Thereafter in section \ref{Parcing} a parcing formula will be given for quasismooth projectively Groenstein varieties with orbifold loci of dimension $\leq 1$, and the fomula still keeps the integrality and symmetric properties as in \cite{BRZ}.

\smallskip
\textit{Acknowledgement} I am gratefull to my suppervisor Miles Reid for guidance and help through my Ph.D study. 
\section{Notations}\label{Notations}
We work over complex field $\CC$. Our definitions and notations about subvarieties of weighted projective space can mostly be found in \cite{Fl}.  We will transfer some the notations to the stack setting, and we refer to \cite{DM} and the appendix of \cite{Vi} for the definitions and notations about stacks. When it comes to parcing Hilbert series, we follow the same convention as in \cite{BRZ}.

A weighted projective space $\PP(a_0,\ldots,a_n)$ is given by the $\CC^\times$ quotient $\CC^{n+1}\setminus \{0\}/\CC^\times$, where the $\CC^\times$ quotient is given by 
\[
(x_0,\ldots,x_n)\rightarrow (\lambda^{a_0}x_0,\ldots, \lambda^{a_n}x_n), \,\, \text{for all }\lambda \in \CC^\times,
\]
where $x_0,\ldots,x_n$ are the coordinates of $\CC^{n+1}$. On the other hand, it can be also given by $\Proj k[x_0,\ldots, x_n]$ with $\wt x_i=a_i$. 

Let $X$ be a porjective variety. Given a polarization $D$ ($D$ is an ample, $\QQ$ Cartier divisor) on $X$, we have an embedding of $X=\Proj R$, where $R=\oplus_{i\geq 0} H^0(X, \Oh(iD))$, in weighted projective spaces.  $(X, D)$ is called quasismooth if the affine cone $\Spec R $ is smooth outside the origin. When $(X,D)$ is quasismooth, the singularities of $X$ all comes from the global $\CC^\times$ quotient. Due to this reason, we refer them as orbifold loci.
\begin{defn}
Let $\pi: C(X) \rightarrow X$ be the $\CC^\times$ quotient map. A geometric point $P$ in $X$ is said to be an \textit{orbifold point} (or an \textit{orbipoint}) \textit{of type $\frac{1}{r}(b_1,\dots, b_n)$ }if the isotropy group (or the stabilizer group) of $\pi^{-1}(P)$ is $\mu_r$, and $\pi^{-1}(P)$ has local parameters $y_1,\dots, y_n$, such that $\CC^\times$ acts on $y_1,\dots, y_n$ with weights $b_1, \dots, b_n$ respectively. We can also talk about higher dimensional orbi\-fold loci. A curve $\sC\subset X$ is said to be an\textit{ orbicurve of type $\frac{1}{s}(c_1,\dots,c_{n-1})$} if every generic point on $\sC$ is of type $\frac{1}{s}(0,c_1,\dots, c_{n-1})$ (possibly after change of coordinates). Similarly for higher dimensional orbifold loci. 
\end{defn} 

\begin{defn}
A projective variety $(X, D)$ is called projectivley Gorenstein if the ring $R=\oplus_{i\geq 0} H^0(X, \Oh(iD))$ is a Gorenstein ring.  
\end{defn}
Given a projectively Gorenstein variety $(X,D)$, one has $\omega_X=\Oh(k_XD)$ for some $k_X\in \ZZ$, called the canonical weight of $(X,D)$. The Hilbert series $P(t)=\sum_{i\geq 0} h^0 (X,\Oh(iD))t^i$ , associated to a projectively Gorenstein varieity $(X,D)$  has the property of Gorestein symmetry, i.e. it satisfies $P(1/t)=(-1)^{n+1}t^{-k_X}P(t)$ (see Lemma 1.1, \cite{BRZ} for more details).

We define weighted projective stack $\underline{\PP}(a_0,a_1,\cdots,a_n)$ as the quotient stack $[\CC^{n+1}\setminus \{0\}/\CC^*]$, where $\CC^*$ action is the same as for weighted projective space $\PP(a_0,a_1,\cdots,a_n)$. We refer $\Stac R$ as the quotient stack $[\Spec R/\CC^*]$. Then on $\Stac R$ we can define similarly quasismoothness and the type of its orbifold loci as above.

\section{Riemann--Roch theorem}\label{Riemann-Roch}

In \cite{To}, Theorem $4.10$, To{\"e}n gave a Riemann--Roch formula for sheaves on smooth Deligne--Mumford stacks. In this section we are going to first recall the ideas of the proof of this general Riemann--Roch theorem, and then translate it to our case where the stacks concerned are quasismooth substacks of weighted projective stacks.  

\subsection{Idea of the Riemann--Roch formula}

Here we will go through the argument working with Deligne--Mumford quotient stacks for simplicity and also because we are mainly concerned with this type of stack. 

Given a Deligne--Mumford stack $\sX$, let $\text{Vect}(\sX)$ (repectively, $\text{Coh}(\sX)$) be the category of vector bundles (respectively, coherent sheaves) on $\sX$. In \cite{To}, To{\"e}n uses Quillen's higher $K$-theory \cite{Qui}, which defines $K_*(\sX)$ to be the homotopy groups of the classifying space $BQ\text{Vect}$ and $G_*(\sX)$ to be the homotopy goups of $BQ\text{Coh}$, see \cite{Qui} for details. Theroem $1$ of \cite{Qui} says $K_0(\text{Vect}(\sX))$ is canonically isomorphic to the Grothendieck group $K_0$, i.e. the free group genetated by vector bundles on $\sX$ modulo the relation induced by exact sequences. For orbifolds, we know that every coherent sheaf is a quotient of a vector bundle and therefore the natural morphism $K_0(\sX)\rightarrow G_0(\sX)$ is an isomorphism. 

Next, we need to set up the link between vector bundles on $\sX$ and vector bundles on its inertia stack $I_{\sX}$.  Theorem $3.15$ in \cite{To} defines a map 
\[
\phi\colon K_0(\sX)\rightarrow K_0(I_{\sX})\otimes \Lambda,
\]
where $\Lambda=\QQ[\mu_{\infty}]$ and $\mu_{\infty}$ is the group of all the roots of unity. This map is the composition of two maps. The first is $ \pi^*\colon K(\sX)\rightarrow K(I_{\sX})$, where $\pi$ is the natural map $\pi:I_{\sX}\rar \sX$. Recall that a vector bundle $\sV$ on $I_\sX$ is given by the following data: 
\begin{itemize}
	\item To every section $s: U \rightarrow \sX$ and every automorphism $\alpha \in \text{Aut}(s)$, where $U \in \text{Sch/S}$, one associates a vector bundle $\sV_{s,\alpha}$ over $U$.
	\item For every pair $(s,\alpha)$ in $I_{\sX}(U)$ and $(s^{\prime},\alpha^{\prime})$ in $I_{\sX}(V)$, every morphism $f:V\rightarrow U$ of $S$-schemes, and every isomorphism $H:f^*(s,\alpha)\cong (s^\prime, \alpha^\prime)$, there is an isomorphism of vector bundles: 
	\[
	\varphi_{f,H}:f^*\sV_{s,\alpha}\cong \sV_{s^\prime, \alpha^\prime}.
	\]
	\item For all pair of morphisms of $S$-schemes
	\[
	W\xrightarrow{g} V
	 \xrightarrow{f} U,
	\]
  all objects $(s,\alpha)$ in $I_{\sX}(U)$, $(s^{'}, \alpha^{'})$ in $I_{\sX}(V)$, $(s^{''}, \alpha^{''})$ in $I_{\sX}(W)$ and all isomorphisms $H_1:f^*(s,\alpha)\cong (s^{'}, \alpha^{'})$  and $H_2: g^* (s^{'}, \alpha^{'}) \cong  (s^{''}, \alpha^{''})$, there is an equality:
  \[
  g^*\varphi_{f,H_1}\circ \varphi_{g,H_2} \cong \varphi_{f\circ g, g^*H_1\circ H_2}.
  \]
	
\end{itemize}
Then $\pi^*:K(\sX)\rightarrow K(I_{\sX})$ can be given as follows: for any vector bundle $\sV$ on $\sX$,  $(\pi^*\sV)_{s,\alpha}$ on all pairs $(s,\alpha)$, with $s:U \rightarrow \sX$ and $\alpha \in \text{Aut}(s)$, are all given by the sheaf $\sV_U$ of $\sV$ on the section $s:U\rightarrow \sX$.

The second map of the composition $\phi$ is the map $\text{dec}\colon K(I_{\sX}) \rightarrow K(I_{\sX})\otimes \Lambda$ which decomposes sheaves into their eigensheaves. In fact, for all objects $(s,\alpha)\in I_{\sX}(U)$ and automorphisms $\alpha$ of $s$ in $\sX(U)$, $\alpha$ defines an isomorphism $H \colon (s,\alpha)\rar (s,\alpha)$ in $I_{\sX}(U)$. Therefore by the above description, a vector bundle $\sV_{(s,\alpha)}$ on $U$ comes naturally with an action of the cyclic group $\langle \alpha \rangle$. Since $\alpha$ is of finite order $r$, the action can be diagonalized canonically as $\sV_{(s,\alpha)}\cong \sV_{(s,\alpha)}^{(\ep)}\bigoplus W_{(s,\alpha)}$, where $\alpha$ acts on $\sV_{(s,\alpha)}^{(\ep)}$ by multiplication of $\ep$, and $\ep$ is in the $r$-th roots of unity. In this way, one can define a subbundle $\sV^{(\ep)}$ of $\sV$ on $I_\sX$. The map $\text{dec}$ sends every vector bundle $\sV$ to the sum of eigen subbundles $\bigoplus_{\ep \in \mu_{\infty}}\ep \sV^{(\ep)}$. 

Combining these two maps $\pi^*$ and $\text{dec}$, we get $\phi=\text{dec}\circ \pi^*: K_0(\sX) \rightarrow K_0(I_{\sX})\otimes \Lambda$ which sets up the link between the $K_0$-theory of the stack and $K_0$-theory of its inertia stack. 
\smallskip

These maps can be given explicitly for quotient stacks. Recall that for a quotient stack $[X/G]$, a sheaf on $[X/G]$ is equivalent to a $G$-equivariant sheaf on $X$, and $I_{[X/G]}$ is isomorphic to $\bigsqcup_{g\in G}[X^g/G]$, where $X^g$ is the fixed locus of $g$ for every $g \in G$. Given an $G$-equivariant vector bundle $\sV$ on $X$, then $\sV$ restricted to the fixed locus $X^g$ is still $G$-equivariant on $X^g$ for any $g$ since $X^g$ is invariant under the $G$ action. Therefore $\sV$ is mapped to a sheaf on $I_{[X/G]}$ by restricting to each component of the inertia stack and one can check that this is the same as $\pi^*\sV$. Given an equivariant vector bundle $\sV$ on $X^g$ for some $g\in G$, $g$ acts on $X^g$ trivially and thus $g$ acts on the fibers of the vector bundle. Thus $\sV$ can be decomposed into eigensheaves $\sV=\bigoplus_{\ep \in \mu_r} \sV^{(\ep)}$, where $g$ acts on the subsheaf $\sV^{(\ep)}$ through multiplication by $\ep$. In this case, $\text{dec}$ sends each $\sV$ to the direct sum $\bigoplus_{\ep \in \mu_r} \ep \sV^{(\ep)}\in K([X^g/G])\bigotimes \Lambda$.

\smallskip

One more concept we need to set up is the conormal bundle of the inertia stack $\sN^*$. In the case of a quotient stack $[X/G]$, this notion is straightforward since each component of the inertia stack $\bigsqcup_{g\in G}[X^g/G]$ is naturally embedded in the originally stack $[X/G]$ and therefore the conormal bundle of each component in $[X/G]$ is well defined. In fact, the tangent sheaf $\sT_{[X/G]}$ comes from a equivariant sheaf of $X$ and it is naturally equivariant when restricted on $X^g$. The tangent sheaf of $[X^g/G]$ also results from an equivariant sheaf on $X^g$. Therefore the quotient of these two tangent sheaves is still equivariant on $X^g$, which defines the normal bundle of $[X^g/G]$ in $[X/G]$. In this way we obtain the normal bundle of the inertia stack $I_{[X/G]}$ in $[X/G]$.

Now let $\alpha_{\sX}= \text{dec }(\lambda_{-1}(\sN^{*}))$, where $\lambda_{-1}(\sN^{*})= \sum (-1)^i \wedge^{i} \sN^*$ as in \cite{FL}. Then Riemann--Roch can be obtained by combining the following two diagrams. Let $\sX$ and $\sY$ be two smooth stacks. For every proper morphism $f:\sX \rar \sY$ the following diagram given in Lemma $4.11$ in \cite{To} commutes:
\begin{displaymath}
    \xymatrix{
       K_0(\sX)\ar[rr]^{\alpha_\sX^{-1}\phi} \ar[d]_{f_*} &&K_0(I_{\sX})\bigotimes \Lambda  \ar[d]^{If_*} \\
       K_0(\sY) \ar[rr]_{\alpha_\sY^{-1}\phi}               && K_0(I_{\sY})\bigotimes \Lambda  }
\end{displaymath}
where $f_*$ is given by $\sum_i(-1)^i R^if_*(-)$ and $If$ is induced by $f$. Another commutative diagram given in Lemma $4.12$ in \cite{To} is the following:
\begin{displaymath}
    \xymatrix{
       K_0(I_{\sX})\otimes \Lambda \ar[rr]^{\ch(-)\td_{I_{\sX}}} \ar[d]_{If_*} &&A(I_{\sX})\otimes \Lambda  \ar[d]^{If_*} \\
       K_0(I_{\sY}) \otimes \Lambda  \ar[rr]_{\ch(-)\td_{I_{\sY}}}            && A(I_{\sY}) \otimes \Lambda }
\end{displaymath}
where $\ch$ and $\td$ are the Chern character and the Todd character which can be defined in the usual way. Here one can take $A(I_{\sX})$ or $A(I_{\sY})$ to be the rational Chow group of $I_{\sX}$ or $I_{\sY}$ defined in \cite{Vi}, Definition $3.4$. Combining these two commutative diagrams, we arrive at the Grothendieck Riemann--Roch theorem obtained by To{\"e}n.
\begin{theorem}(B.To{\"e}n)\label{RR}
Let $\sX$ and $\sY$ be smooth stacks. Define the representation Todd class $\td^{\text{rep}}_\sX$ to be $\ch(\alpha_\sX^{-1})\td_{I\sX}$. Then for any $\sF \in K_0(\sX)$ and any proper morphism $f:\sX\rar \sY$, one has:
\[
If_*(\ch(\phi(\sF))\td_\sX^{\text{rep}})=\ch (\phi(f_*(\sF))\td_\sY^{\text{rep}}.
\]
\end{theorem}

\subsection{Riemann--Roch formula for quasismooth projective stacks}

To write down the Riemann--Roch formula for quasismooth stacks, we need to introduce some more notation (see \cite{N}).

Let $\sX$ be a quasismooth projective substack $\Stac R$ inside $\underline{\PP}(a_0,\dots,a_n)$, where $R=k[x_0,\dots,x_n]/J$ with $J$ a weighted homogeneous ideal. Let $I_{\underline{\PP}}$ (resp.\ $I_{\sX}$) be the inertia stack of $\underline{\PP}$ (resp.\ $\sX$). Then there is a natural embedding $I_{\sX} \hookrightarrow I_{\underline{\PP}}$, which is given in each component of $I_{\underline{\PP}}$, say $\underline{\PP}(a_{i_0},\dots,a_{i_m})$, by the substack $\sY=\Stac R^\prime$, where $R^\prime=k[x_{i_0},\dots,x_{i_m}]/J\cap k[x_{i_0},\dots, x_{i_m}]$. 

Let $S=\{\text{all subsets of } \{a_0,\dots, a_n\}\}$. The subset $S_0$ of $S$ is defined as follows:
$$
S_0 = \left\{  \{a_{i_0},\dots, a_{i_m}\} \in S \left|
\begin{array}{l}
\nexists {a_{j_0},\dots, a_{j_l}}\, s.t. \\
 \mathrm{gcd} (a_{j_0},\dots, a_{j_l}) = \mathrm{gcd} (a_{i_0},\dots, a_{j_m})\, \mathrm{and} \\
\{a_{i_0},\dots,a_{i_m} \} \subset \{a_{i_0},\dots, a_{i_l} \}
\end{array} \right. \right \}.
$$
In other words, it contains the subsets of $\{a_0,\dots, a_n\}$ which are the largest among these who have the same greatest common divisors. For instance, let $S=\{1,3,4,6\}$ then $S_0=\{\{1,3,4,6\}, \{4,6\}, \{3,6\},\{4\},\{6\}\}$. Moreover, for each of the subsets $s=\{a_{i_0},\dots, a_{i_m}\} \in S_0$ with $r=\text{ gcd }(a_{i_0},\dots,a_{i_m})$, we associate to it a set $\tau_s$, which is defined by
$$
\tau_s=\left\{  \ep \in \mu_r \left|
\begin{array}{l} \ep \notin \mu_q, \text{ if there exists } \{a_{j_0},\dots,a_{j_l}\} \in S_0
\,\, s.t. \\
q=\text{ gcd }(a_{j_0},\dots,a_{j_l}) \text{ and } q|r
\end{array}\right.\right\}.
$$
Take the above example. To $s=\{6\}\in S_0$, we associate the set $\{\ep \in \mu_6\,|\, \ep^2 \neq 1 \text{ and }\ep^3 \neq 1\}$. Using these notation, the inertia stack $I_{\underline{\PP}}$ is given by $\sqcup_{s\in S_0} (\underline{\PP}(s) \times \tau_s)$, where $\underline{\PP}(s)=\underline{\PP}(a_{i_1},\dots,a_{i_m})$ and $x_{i_j} \in s$ with weight $a_{i_j}$. If we let $\sY_s$ be the substack of $\underline{\PP}(s)$ defined by the ideal $J$, then the inertia stack $I_{\sX}$ of $\sX$ is given by $\sqcup_{s\in S_0}(\sY_s  \times \tau_s)$. 

Using above notations, we can state the Riemann--Roch formula for quasismooth stacks. 
\begin{prop}\label{RRforquasismooth}
Let $\sX$ be a quasismooth substack in a weighted projective stack $\underline{\PP}(a_0,\dots,a_n)$, and let $\sV$ be a vector bundle on $\sX$. Using the above notation, one has 
\[
\chi(\sV)=\sum_{s \in S_0} \sum_{\ep \in \tau_s} \bigl[\frac{\ch(\phi(\sV))\td \sY_s}{\ch(\lambda_{-1}(\mathrm{dec }\,(\sN_s^*)))}\bigr]_{\mathrm{dim }\,\sY_s },
\]
where $\sN_s^*$ is the conormal bundle of $\sY_s$ inside $\sX$ and $[-]_{\mathrm{dim}\, \sY_s}$ represents the codimension $\mathrm{dim }\, \sY_s$ part in the Chow group.  In particular, when $\sV=\Oh(d)$, then
 
\[
\chi(\Oh(d))=\sum_{s\in S_0} \sum_{\ep \in \tau_s} \bigl[\frac{\ep^d \ch(\Oh(d))\td \sY_s}{\ch(\lambda_{-1}(\mathrm{dec}\,(\sN_s^*)))}\bigr]_{\mathrm{dim}\, \sY_s}. 
\]

\end{prop}
\begin{pf}
In Theorem \ref{RR}, if we take $\sY$ to be a point, we will get the Hizebruch--Riemann--Roch formula for a vector bundle $\sV \in K_0(X)$. In the first diagram above Theorem $\ref{RR}$, the map $\alpha_{\sX}^{-1}\phi$ sends $\sV$ to a direct sum of sheaves on $I_{\sX}$, and in the second diagram, we can calculate $\ch$ and $\td$ componentwise on $I_{\sX}$. Then we obtain the Riemann--Roch formula for vector bundles on $\sX$. In particular, if $\sV=\Oh(d)$, then for each $\sY_s$ and each element $\ep \in \tau_s$, one has $\phi(\sV)=\ep^d\sV|_{\sY_s}$.   $\square$ 
\end{pf}

\smallskip

\begin{rem}
Let $\sY_s$ be one of the components of the inertia stack of $I_{\sX}$ and $\tau_s$ the set associated to it. Suppose the normal bundle $\sN_{s}$ of rank $r$ in $\sX$ of $\sY_s$ can be decomposed into the direct sum $\bigoplus_{i=1}^{l} \sN_i$ under the group $<\ep>$ action for each $\ep \in \mu_r$, and each $\sN_i$ has eigenvalue $\ep^{-a_i}$, then the denominator of the formula in the proposition can be written as
\begin{eqnarray*}
&&\ch(\lambda_{-1}(\text{dec} (\sN^*)))=\ch(\lambda_{-1}(\bigoplus_{i=1}^{l}\ep^{-a_i}\sN_{i}^*))\\
&&=\ch(\prod_{i=1}^{l}(1-\ep^{-a_i}\sN_i^*))=\prod_{i=1}^{l}(1-\ep^{-a_i}e^{-v_i}),
\end{eqnarray*}
where $v_i$ is the first Chern class of $\sN_i$. Moreover, we can express the inverse
\begin{eqnarray*}
&&\frac{1}{(1-\ep^{-a_i}e^{-v_i})}=\frac{1}{1-\ep^{-a_i}}-\frac{\ep^{-a_i}}{(1-\ep^{-a_i})^2}v_i+\\
&&(\frac{\ep^{-a_i}}{(1-\ep^{-a_i})^3}-\frac{\ep^{-a_i}}{2(1-\ep^{-a_i})^2})v_i^2+\text{higher order terms }.
\end{eqnarray*}
This expression is very useful, as we will see in the concrete cases below. 
\end{rem}

\bigskip

Using the formula in Proposition \ref{RRforquasismooth} and the above remark, for quasismooth stacks with concrete orbifold loci one can express this formula in terms of Dedekind sums. Given a quasismooth stack $\sX$ of dimension $n$ with only isolated orbipoints $\sB=\{P \text{ of type }\frac{1}{r}(b_1,\dots,b_n)\}$, the inertia stack of $\sX$ can be written as $I_{\sX}=\sX \sqcup_s (\sY_s\times \tau_s) $ where $\sY_s$ are all of dimension $0$, corresponding to the orbifold points, and $\tau_s$ is $\mu_r\setminus \{1\}$, determined by the orbifold type of $\sY_s$. Then the formula in Proposition \ref{RRforquasismooth} can be written as:
\begin{cor}\label{RRisolatedsingularities}
Given $\sX$ as above, the Riemann--Roch formula for $\Oh_{\sX}(d)$ is given by 
\begin{equation}\label{RRforisolatedorbipoints}
\chi(\Oh_{\sX}(d))=[\ch(\Oh_{\sX}(d))\td_{\sX}]_n+\sum_{P\in \sB}\frac{1}{r}\sum_{\ep \in \mu_r,\ep \neq 1}\frac{\ep^d}{\prod_i(1-\ep^{-b_i})}.   
\end{equation}
\end{cor}
\begin{pf}
In this case, the only components for the inertia stack are the stack itself and the orbipoints. Each of the orbipoints of type $\frac{1}{r}(b_1,\cdots,b_n)$ is associated with $r-1$ components of the inertia stack, namely $\sqcup_{\ep \in \mu_r,\ep \neq 1}[C(P)/\CC^*]\times \ep$, where $C(P)$ is the orbit of $P$. Now consider one of the components $[C(P)/\CC^*]\times \ep$ corresponding to a singular point of type $\frac{1}{r}(b_1,\dots, b_n)$ with normal bundle $\sN$. Then $\ch (\lambda_{-1}(\text{dec }(\sN^*)))$ is equal to $\ch(\lambda_{-1}(\text{dec }(\bigoplus \sN_{i}^*)))=\prod_i(1-\ep^{-b_i}e^{-v_i})$, where $v_i$ is the first Chern class of $\sN_i$. Since each of the $\sY_s$ is of dimension $0$, we have that $\ch(\phi(\Oh(d)))=\ep^d$ and $\td_{\sY_s}=1$, and for each component $[C(P)/\CC^*]\times \ep$, we have 
\[
\bigl[\frac{\ep^d \ch(\Oh(d))\td_{\sY_s}}{\ch(\lambda_{-1}(\text{dec}(\sN_s^*)))}\bigr]_{0}=\frac{1}{r}\frac{\ep^d}{\prod_i(1-\ep^{-b_i})},
\]
where $\frac{1}{r}$ is the degree of the point. Summing over all the components we get the formula.  $\square$
\end{pf}

\begin{rem}\label{rem}
Note for $d=0$, one obtains
\[
\chi(\Oh_{\sX})=\td_n+\sum_{P\in \sB} \frac{1}{r}\sum_{\ep \in \mu_r,\ep \neq 1} \frac{1}{(1-\ep^{-b_1})\cdots(1-\ep^{-b_n})}, 
\]
where $\td_n$ represents the top Todd class of $\sX$. Thus replacing the $\td_n$ via the above equality (\ref{RRforisolatedorbipoints}) gives the same formula as in \cite{R}.
\end{rem}

\smallskip
Now suppose $\sX$ has orbifold loci of dimension $\leq 1$, and its orbifold loci are 
\begin{itemize}
 \item the set of all orbicurves $\sB_C=\{\text{ orbicurves of type } \frac{1}{r}(a_1,\dots, a_{n-1})\}$, and
 \item the set of all orbipoints $\sB_P=\{\text{orbipoints of type }\frac{1}{s}(b_1,\dots,b_n)\}$.
\end{itemize} 
In this case, we have $I_{\sX}=\sX\sqcup_{\sB_C} (\sqcup_{\ep \in \mu_r, \ep 
\neq 1} \sC \times \ep)\sqcup_{\sB_P} (\sqcup_{\ep \in \mu_s, \ep^{b_i}\neq 1} P\times \ep )$ and the Riemann--Roch formula is given by 
\begin{cor}\label{RRformula}
Given such an $\sX$ with only orbifold loci of dimension $\leq 1$, one has 
\[
\chi(\Oh_{\sX}(d))=[\ch(\Oh_{\sX}(d))\td_{\sX}]_n +\sum_{P\in \sB_P} M_P+ \sum_{\sC\in \sB_C} M_{\sC},
\]
where $M_P$ for a point $P$ of type $\frac{1}{s}(b_1,\dots, b_n)$ is given by
\[
\frac{1}{s}\sum_{\ep \in \mu_r,\ep^{-b_i} \neq 1}\frac{\ep^d}{\prod_i(1-\ep^{-b_i})},
\]
while $M_{\sC}$ for a curve $\sC$ of type $\frac{1}{r}(a_1,\dots,a_{n-1})$ is given by
\begin{eqnarray*}
&&\frac{1}{r}\sum_{\ep \in \mu_r,\ep \neq 1} \frac{\ep^d}{\prod(1-\ep^{-a_i})} d \deg H|_{\sC}-\frac{1}{2r}\sum_{\ep \in \mu_r,\ep \neq 1}  \frac{\ep^d}{\prod(1-\ep^{-a_i})}\deg K_{\sC}\\
&&-\sum_{i=1}^{n-1}\frac{1}{r}\sum_{\ep \in \mu_r,\ep \neq 1} \frac{\ep^{d-a_i}}{(1-\ep^{-a_i})^2 \prod_{j \neq i} (1-\ep^{-a_j})} \deg \gamma_i,
\end{eqnarray*}
where $H $ is the divisor (possibly $\QQ$-divisor) corresponding to the sheaf $\Oh_{\sX}(1)$, and $\gamma_i$'s are the first chern classes of $\sN_i$ in the decomposition of the normal bundle $\sN=\bigoplus_{i}\sN_i$.
\end{cor}
\begin{pf}
As in the proof of Corollary \ref{RRisolatedsingularities} we obtain the part coming from orbifold points $M_P$. An orbicurve of type $\frac{1}{r}(a_1,\dots, a_{n-1})$ will give rise to $r-1$ components in the inertia stack of $\sX$, namely, $\sqcup_{\ep \in \mu_r,\ep \neq 1} \sC\times \ep$. We also know that the normal bundle of each component can be decomposed into $\bigoplus \sN_i$ with $\sN_i$ be in the $\ep^{a_i}$ eigenspace. Suppose $c_1(\sN_i)=\gamma_i$. Then in the formula, for the component $\sC\times \ep $ we will have 
\begin{eqnarray*}
&&\bigl[\frac{\ep^d\ch(\Oh(d))\td_{\sC}}{\ch(\lambda_{-1}(\sN^*))}\bigr]_1\\
&=&[(1+d H|_{\sC})(1+\frac{1}{2}c_1(\sT_C))\prod_{i=1}^{n-1}(\frac{1}{1-\ep^{-a_i}}-\frac{\ep^{-a_i}}{(1-\ep^{-a_i})^2}\gamma_i)]_1\\
&=& \frac{\ep^d(dH|_{\sC}+\frac{1}{2}c_1(\sT_C))}{\prod_{i=1}^{n-1}(1-\ep^{-a_i})}-\sum_{i=1}^{n-1}\frac{\ep^{d-a_i}}{(1-\ep^{-a_i})^2\prod_{j\neq i}(1-\ep^{-a_j})}\deg \gamma_i,
\end{eqnarray*}
where $H$ is the $\QQ$-divisor corresponding to $\Oh(1)$ and $\sT_{\sC}$ is the tangent sheaf of $\sC$. Summing these over the $r-1$ components in the inertia stack, we get the above formula.  $\square$
\end{pf} 
\begin{rem}
In the above formula, by abuse of notation, we write $\deg H|_{\sC}$ for the number given by the intersection number of $rH$ with $\sC$, because in this way the coefficients can be given in the form of Dedekind sums as in Section \ref{app}. Similarly for $\deg K_{\sC}$, the $\deg \, K_{\sC}$ here is given by $r$ times degree of the divisor $K_{\sC}$, where $K_{\sC}$ is the canonical divisor of $\sC$ as a stack. For example, $\sC=\PP(2,4)$ has $\deg K_\sC=2\times (-\frac{6}{8})=-\frac 32$. We will also use the same convention in the following. 
\end{rem}

Of course, we can continue to write out the formula for quasismooth stacks with orbifold loci of dimension $\geq 2$ in the same way, but we will omit her.

\subsection{Riemann--Roch on the moduli space}
In the last section, we obtained the Riemann-Roch formula for line bundles $\Oh_{\sX}(d)$ on $\sX=\Stac R$. Now we want to deduce the Riemann-Roch formula for $\Oh_{X}(d)$ on its moduli space $X= \Pro R$. For this we just need to set up the link between $\sX$ and $X$.

Let $\pi:\sX\rightarrow X$ be the map induced by the quotient map $\hat{\pi}:\mathrm{Spec} R\setminus\{0\} \rightarrow X$. Then $\pi$ is the natural map from $\sX$ to $X$ inducing a bijection between the geometric points of $\sX$ and $X$. Recall that we define $\Oh_{\sX}(d)$ to be the line bundle descended from an equivariant line bundle on the affine cone, but we can also define it on an \'{e}tale cover of $\sX$, in which case we can see clearly that $\pi_*(\Oh_{\sX}(d))=\Oh_X(d)$. Calculating the \v{C}ech cohomology on $\sX$ and $X$ gives us $H^i(\sX,\Oh_{\sX}(d))=H^i(X,\Oh_{X}(d))$ for all $i$, and therefore $\chi(\sX,\Oh_{\sX}(d))=\chi(X,\Oh_{X}(d))$.    

In this way, we can transfer the formula for $\chi(\Oh_{\sX}(d))$ to the coarse moduli space $X$ to get a formula for $\Oh_{X}(d)$. Recall that the formula for $\chi(\Oh_{\sX}(d))$ is given by a sum over all the components of the inertia stack $I_{\sX}$, which implies that the formula on $X$ will sum over all the singular strata of $X$.  We also know that the morphism $\pi_*:A(\sX)\otimes \QQ \rightarrow A(X)\otimes \QQ$ between Chow groups given in \cite{Vi} is an isomorphism. For an integral closed substack of $\sY$, the map $\pi_*$ sends $[\sY]$ to $[\frac{1}{g_{\sY}}\pi(\sY)]$, where $g_{\sY}$ is the order of the generic stabilizer group of $\sY$. 
\bigskip

When the quasismooth stack $\sX=\Stac R$ has only codimension $\geq 2$ orbifold loci, then the coarse moduli space given by $X=\Pro R$ has cyclic quotient singularities in one to one correspondence with the orbifold loci on $\Stac R$. Take the case when there are only curve singularities as an example. 

\begin{prop}\label{RRwithcurvesingularities}
Let $X$ be a quasismooth variety of dimension $\geq 3$ in weighted projective space $\PP(a_0,\dots,a_n)$. Let $\sB=\{$C$ \text{ singular of type }\\
\frac{1}{r}(a_1, \dots, a_{n-1})\}$ be all the singular loci on $X$. Then 
\[
\chi(\Oh(d))=[\ch(\Oh(d))\td_X]_n+ \sum_{C\in \sB} M_C
\]
where $M_C$ is given by
\begin{eqnarray*}
&&\frac{1}{r}\sum_{\ep \in \mu_r, \ep \neq 1} \frac{\ep^d}{\prod(1-\ep^{-a_i})} d \deg H|_{C}-\frac{1}{2r}\sum_{\ep \in \mu_r, \ep \neq 1}  \frac{\ep^d-1}{\prod(1-\ep^{-a_i})}\deg K_{X}|_{C}\\
&&-\sum_{i=1}^{n-1}\frac{1}{2r} \sum_{\ep \in \mu_r, \ep \neq 1} \frac{(\ep^{d}-1)(1+\ep^{-a_i})}{(1-\ep^{-a_i})^2 \prod_{j \neq i} (1-\ep^{-a_j})} \deg \gamma_i.
\end{eqnarray*}
where $H$ is the Weil divisor associated to $\Oh_X(1)$, and the $\gamma_i$ are the first Chern classes of the orbibundle $\sN_i$, and $\bigoplus_{i=1}^r\sN_i$ is the decomposition of the normal bundle $\sN$ of $C$ in $X$.  
\end{prop}
\begin{pf}
Here we just need to point out that intersection number $\deg K_{X}|_C$ is defined as follows: Let $\hat{C}=\hat{\pi}^{-1}(C)$. Then $\sC=[\hat{C}/\CC^*]$ is a substack of $\sX$, which maps to $C$ by $\pi$. Since $\pi^* K_{X}=K_{\sX}$, by projection formula, we have 
\[
K_{\sX}\cdot \sC=\pi_*(K_{\sX}\cdot \sC)=\pi_*(\pi^*K_X\cdot \sC)=K_X\cdot \frac{1}{r}C.
\]
Similarly for $H|_C$. $\square$
\end{pf}
\begin{rem}
Here the definition of $\deg H|_C$ coincides with the one given in \cite{BS}. Therefore as a special case we can recover the formula in \cite{BS}.
\end{rem}

\subsection{Calculating Dedekind sums}\label{app}
Before going any further, we would like to study the Dedekind sums appeared in the formulas so that we will be able to characterize and calculate them. Here by Dedekind sum, we mean a sum of the form:
\begin{eqnarray*}
\sigma_i(\frac{1}{r}(a_1,\dots,a_n))&=&\frac{1}{r}\sum_{\ep \in \mu_r, \ep^{a_i}\neq 1}\frac{\ep^i}{(1-\ep^{-a_i})\dots(1-\ep^{-a_n})}\\
&=&\frac{1}{r}\sum_{\ep \in \mu_r, \ep^{a_i}\neq 1}\frac{\ep^{-i}}{(1-\ep^{a_i})\cdots(1-\ep^{a_n})},
\end{eqnarray*}
where $(a_1,\dots,a_n)$ is a sequence of positive integers such that $a_i$ mod $r \neq 0$ for all $i$. Such sums are closely related to traditional Dedekind sums, thus we still refer it as the $i$th Dedekind sum, denoted by $\sigma_i(\frac{1}{r}(a_1,\dots,a_n))$ or simply $\sigma_i$.  We write $\delta_i$ for $\sigma_i-\sigma_0$, that is,
\[
\delta_i=\frac{1}{r}\sum_{\ep \in \mu_r, \ep^{a_i}\neq 1}\frac{\ep^{-i}-1}{(1-\ep^{a_i})\cdots(1-\ep^{a_n})}.
\]
When $n=1$ and $(a,r)=1$, there is a compact expression for $\delta_i(\frac{1}{r}(a))$.
\begin{lem} \label{line}
When $(a,r)=1$,
\[
\delta_i(\frac{1}{r}(a))=\sigma_{i}(\frac{1}{r}(a))-\sigma_0(\frac{1}{r}(a))=\frac{1}{r}\sum_{\ep\in \mu_r,\ep \neq 1} \frac{\ep^{-i}-1}{1-\ep^{a}}=-\frac{\overline{bi}}{r},
\]
where $b$ is the inverse of $a$ modulo $r$, i.e., $ab=1$ mod $r$. In particular, this gives
\[
\sigma_0(\frac{1}{r}(a))= \frac{r-1}{2r}.
\]
\end{lem}
\begin{pf}
Let $ab=1$ mod $r$, then $(\ep^{a})^{\overline{bi}}=\ep^i$, where $\overline{bi}$ represents the smallest nonnegative residue of $bi$ modulo $r$ (similarly in what follows). Thus
\begin{eqnarray*}
\ep^{r-i}-1&=&(\ep^{a})^{r-\overline{bi}}-1\\
           &=&((\ep^{a})^{r-\overline{bi}-1}+\cdots+1)(\ep^a-1).
\end{eqnarray*}
Note that $\sum_{\ep \in \mu_r,\ep\neq 1} \ep^m=-1$ for all $m\neq 0$. Then 
\begin{eqnarray*}
\delta_i=\frac{1}{r}\sum_{\ep\in \mu_r,\ep \neq 1} \frac{\ep^{-i}-1}{1-\ep^{a}}&=& -\frac{1}{r}\sum_{\ep\in \mu_r,\ep \neq 1} ((\ep^{a})^{r-\overline{bi}-1}+\cdots+1)\\
&=&-\frac{1}{r}(\underbrace{(-1+ \cdots+(-1))}_{r-\overline{bi}-1}+r-1)=-\frac{\overline{bi}}{r}.
\end{eqnarray*}
Moreover since $\sum_{i=0}^{r-1}\sigma_i(\frac{1}{r} a)=0$, one has
\[
\sigma_0=\frac{\sum_{i=0}^{r-1} \overline{bi}}{r^2}=\frac{r-1}{2r},
\]
because $b$ is coprime to $r$ and thus $\overline{bi}$ will run over $1,\dots,r-1$ for $0\leq i\leq r-1$.
\end{pf} 
\begin{exa}
Take $r=5$, $a=3$, and one has $b=2$. Thus for $i=1,\dots,4$, the $\delta_i (\frac{1}{5}(a))$ are: $-2/5,-4/5,-1/5,-3/5$.
\end{exa}
To calculate all $\sigma_i$ in general, we have the following proposition to use (see also \cite{B} for a different proof).
\begin{prop}\label{relation}
Given positive integers $r$ and $a_1,\dots,a_n$ such that $a_i$ are not divisible by $r$, let $h=\mathrm{gcd}\,(\prod_{j=1}^{n}(1-t^{a_j}),\frac{1-t^r}{1-t})$. Then $\sum_{i=0}^{r-1} \sigma_{i}t^i$ is the inverse of $\prod_{j=1}^n(1-t^{a_j})$ modulo $\frac{1-t^r}{h(1-t)}$, that is,  
\[
\bigl(\sum_{i=0}^{r-1} \sigma_{i}t^i) \prod_{j=1}^n(1-t^{a_j}\bigr)=1 \,\mathrm{ mod }\, \frac{1-t^r}{h(1-t)}.
\]
\end{prop}
\begin{pf}
Observe that 
\[
\frac{1}{r}\sum_{\ep\in \mu_r,\ep \neq 1}(1+\ep^{-1}\zeta +\cdots +\ep^{-(r-1)}\zeta^{r-1}) \frac{ (1-\zeta^{a_1})\cdots(1-\zeta^{a_n})}{(1-\ep^{a_1})\cdots(1-\ep^{a_n})}=1 
\]
for all $\zeta\in \mu_r, \,\zeta^{a_i}\neq 1,\, \zeta\neq 1$. In fact, when $\ep \neq \zeta$ we have $\sum_{i=0}^{r-1}(\zeta^{-1}\ep)^i=0$ as $\zeta^{-1}\ep$ is still a $r$th roots of unity, and when $\ep =\zeta$ we have $\sum_{i=0}^{r-1}(\zeta^{-1}\ep)^i=r$. Thus we have shown that for all the roots of $\frac{1-t^r}{h(1-t)}$ the left hand side of the equality equals $1$, which is equivalent to:
\[
(\sum_{i=0}^{r-1}\sigma_i t^i)(1-t^{a_1})\cdots(1-t^{a_n})=1 \, \mathrm{mod} \, \frac{1-t^r}{h(1-t)}.
\]
We are done. $\square$
\end{pf}

Using this proposition, we can calculate $\sigma_i(\frac{1}{r}(a_1, \dots, a_n))$ by a computer program. In fact, since $h=\text{ gcd }(\frac{1-t^r}{1-t}, \prod_{j=1}^{n} (1-t^{a_j}))$, by the Euclidean algorithm there exists a unique $\alpha(t)$ of degree $\leq r-\deg h -2$ and $\beta(t) \in \CC[t]$ (in fact, $\alpha(t)$ and $\beta(t)$ are in $\QQ[t]$) such that
\[
\alpha(t) \prod_{i=1}^n (1-t^{a_i})+\beta(t)\frac{1-t^r}{h(1-t)}=1.
\]
This implies that $\alpha(t)$ is also the inverse of $\prod_{i=1}^n(1-t^{a_i})$ modulo $\frac{1-t^r}{h(1-t)}$, and therefore $\alpha(t)=\sum_{i=0}^{r-1} \sigma_i t^i \,\,\mathrm{mod}\,\,\frac{1-t^r}{h(1-t)}$, i.e.,
\[
\sum_{i=0}^{r-1} \sigma_i t^i =\alpha(t) +f(t) \frac{1-t^r}{h(1-t)},
\]  
where $f(t)$ is a polynomial of degree $\deg h$. In particular, $f(t)$ is a constant when $h=1$. If $h\neq 1$, then $f(t)$ will have $\deg h+1$ undetermined coefficients. Thus we need $\deg h+1$ relations among the coefficients of the right hand side to determine $f(t)$ and hence $\sigma_i$. Note that for each $w_i=(a_i,r)\neq 1$ and any $\ep \in \mu_r$, one has $1+\ep^{w_i}+\cdots+\ep^{w_i(r/w_i-1)}=0$. Thus
\[
\sum_{l=0}^{r/w_i-1}\sigma_{w_il+k}=\frac{1}{r}\sum_{\ep\in \mu_r,\ep^{a_i}\neq 1} \frac{(1+\ep^{w_i}+\cdots+\ep^{r-w_i})\ep^k}{(1-\ep^{a_i})\cdots (1-\ep^{a_n})}=0.
\]
for $k=0,1,\dots,w_i-1$. Then for every such $w_i$ there are $w_i-1$ independent relations.  Let $w_{i_j}$ ($j=1,\dots, l$) be all such $w_i$, we know that $\sum_{j=1}^l(w_{i_j}-1)=\deg h$ relations between the $\sigma_i$'s. One more relation comes from the fact that $\sum_{i=0}^{r-1}\sigma_i=0$. Therefore we have in total $\deg h+1$ independent relations among $\sigma_i$, which gives us enough linear equations to determine $f(t)$ and hence $\sigma_i$. This in particular implies $\sigma_i$'s are rational numbers. The following MAGMA program uses above ideas and output $\sigma_0,\dots, \sigma_{r-1}$ if we input $r$ and the sequence $LL=[a_1,\dots,a_n]$.
\begin{pro}\label{contribution}
\begin{verbatim}
function  Contribution(r, LL) 
QQ:=Rationals();
Poly<t>:=PolynomialRing(QQ);
L:=[Integers()|i: i in LL]; n:=#LL;
pi:=&*[(1-t^i):i in L]; A:=Poly!((1-t^r)/(1-t));
G:=GCD(pi, A); dG:=Degree(G);
B:=Poly!(A/G); dB:=Degree(B); 
a,be,c:=XGCD(pi, B); dbe:=Degree(be);
R<[v]>:=PolynomialRing(QQ,dG+2);
va:=Name(R,dG+2); 
bnew:=&+[Coefficient(be,i)*va^i: i in [0..dbe]];
RR:=&+[v[i]*va^(i-1):  i in [1..dG+1]];
Bnew:=&+[Coefficient(B,i)*va^i: i in [0..dB]];
AA:=bnew-RR*Bnew;
S:=[Coefficient(AA,va, 0)] cat [Coefficient(AA, va, r-i): i in [1..
r-1]];
empty:=[];
for a in L do 
dd:=GCD(a,r); tt:=r/dd;
relations:=empty cat [&+[S[dd*l+i]: l in [0..tt-1]]: i in [1..dd]];
empty:=relations;
end for;
Mat:=Matrix(QQ,[[Coefficient(empty[i],v[j],1):j in [1..dG+1]]:i in 
[1..#empty]]);
zero:=[0: i in [1..dG+2]];
V:=-Vector(QQ,[Evaluate(empty[i],zero): i in [1..#empty]]);
MF:=Transpose(Mat); x,y,z:=IsConsistent(MF,V);
yy:=&+[y[i+1]*va^(i):i in [0..dG]];
sigma:=bnew-yy*Bnew;
Sigma:=[QQ!Coefficient(sigma, va, 0)] cat [QQ!Coefficient(sigma, 
va,i): i in [1..r-1]];
return  Sigma;
end function;
\end{verbatim}
\end{pro}
\subsection{Examples}
Now we can do calculations on concrete examples.
\begin{exa}
Consider the subvariety $X_{11}$ of $\PP(1,2,3,5)$, where $X_{11}$ is defined by $f=x_0^{11}+x_1^4x_2+x_1x_2^3+x_0x_3^2$. We can check that it is quasismooth and has $3$ orbipoints $P_1=(0,1,0,0)$, $P_2=(0,0,1,0)$ and $P_3=(0,0,0,1)$ of type $\frac{1}{2}(1,1)$, $\frac{1}{3}(1,2)$, and $\frac{1}{5}(2,3)$ respectively. Hence the formula for the sheaves $\Oh_{X}(d)$ is given by 
\[
\chi(\Oh_{X}(d))=[\ch(\Oh_{X}(d))\td_{X}]_2+ M_{P_1}+M_{P_2}+M_{P_3},
\]
where $M_{P_1}$, $M_{P_2}$ and $M_{P_3}$ are given by Dedekind sums as in (\ref{RRforisolatedorbipoints}), and
\[
[\ch(\Oh_{X}(d))\td_{X}]_2=\td_{2}+\frac{1}{2}dH(dH-K_{X}),
\]
where $H$ is $c_1(\Oh_{X}(1))$. By the exact sequence 
\[
0\rightarrow \sT_{X}\rightarrow \sT_{{\PP}}|_{X}\rightarrow \sN_{X |{\PP}}\rightarrow 0,
\]
we know that $c_t(\sT_{X})c_t(\sN_{X|\PP})=c_t(\sT_{\PP})|_{X}$, and thus we have $c_1(\sT_{X})=-K_{X}=0$ and
\[
c_2(\sT_{X})=c_2(\sT_{\PP}|_{X})-c_1(\sT_{X})c_1(\sN_{X|\PP})=\frac{451}{30}.
\]
Hence $\td_2=\frac{1}{12}(c_1(\sT_{X})^2+c_2(\sT_{X}))=\frac{451}{360}$. Now we use our Program \ref{contribution} to compute $M_{P_i}(d)$ for $0 \leq i\leq 3$.

\begin{verbatim}
>f:=func<d|451/360+1/2*d^2*11/30>;
>MP1:=Contribution(2,[1,1]);
>MP2:=Contribution(3,[1,2]);
>MP3:=Contribution(5,[2,3]);
>[f(d)+MP1[d mod 2 +1]+MP2[d mod 3+1]+MP3[d mod 5+1]: d in
 [1..10]];
[ 1, 2, 3, 4, 6, 8, 10, 13, 16, 20 ]
\end{verbatim}
The last output gives us $\chi(\Oh(d))$ for $1\leq d \leq 10$.
\end{exa}

Next, we give an example with curve orbifold loci and dissident points.
\begin{exa}
Let $X$ be a quasi-smooth Calabi-Yau 3-fold given by $X_{80}\subset \PP^4(3,5,7,25,40)$. It is of degree $2/2625$ and
has an orbifold curve $C_{80} \subset \PP(5,25,40)$ of type $\frac{1}{5}(2,3)$ and a point basket $\sB=\{\frac{1}{3}(1,1,1),\\ \frac{1}{7}(4,5,5),
\frac{1}{25}(3,7,15)\} $, among which the point of type $\frac{1}{25}(3,7,15)$ is a dissident point.
Then according to the Riemann--Roch formula in Corollary \ref{RRformula}, we have several parts in the
formula, which correspond to the connected components of the
associated inertia stack. The first part is given by:
\[
r_1=[\ch(\Oh_X(d)) \td_{X}]_3,
\]
where the Chern character is given by
$\ch(\Oh_{X}(d))=1+dH+d^2H^2/2+d^3H^3/6$. To calculate the
Todd class, we use the exact sequence:
\[
0 \rar \sT_{X} \rar \sT_{\PP}|_{X} \rar \sN_{X|\PP} \rar 0.
\]
Since $X$ is a hypersurface, we have $\sN_{X|\PP}=\Oh_{X}(80)$.
It follows that
\begin{eqnarray*}
c_t(\sT_{X})&=&c_t(\sT_{\PP}|_{X})c_t^{-1}(\sN)\\
&=&(1+3t)(1+5t)(1+7t)(1+35t)(1+40t)(1+80Ht)^{-1}\\
&=&  1+ 2046H^2u^2 - 143960H^3t^3 + \mathrm{higher\,\, order\,\, terms}.
\end{eqnarray*}
That is, $c_1(X)=0,\,c_2(X)=2046H^2,\,c_3(X)=-143960H^3$. 
Thus
\begin{eqnarray*}
r_1&=&[(1+dH+d^2H^2/2+d^3H^3/6)(1+1/2c_1+1/12(c_1+c_2^2)+1/24c_1c_2)]_3\\
&=&1/6d^3H^3+341/2dH^3,
\end{eqnarray*}
where $H^3=\frac{80}{3 \cdot 5 \cdot 7 \cdot 25 \cdot 40}=2/2625 $.

The second part comes from the orbifold curve $C_{80} \subset \PP(5,25,40)$, whose normal bundle is given by $\sN=\Oh_{C}(3)\bigoplus\Oh_{C}(7)$.
Thus the second part $r_2$ is given as follows:
\begin{eqnarray*}
&&\frac{1}{5}\sum_{\ep \in \mu_5,\ep \neq 1} \frac{\ep^d}{(1-\ep^{-2})(1-\ep^{-3})} d \deg H|_{C}-\frac{1}{2\cdot 5}\sum_{\ep \in \mu_5,\ep \neq 1}  \frac{\ep^d}{(1-\ep^{-3})(1-\ep^{-5})}\deg K_{C}\\
&&-\frac{1}{5}\sum_{\ep \in \mu_5,\ep \neq 1} \frac{\ep^{d-3}}{(1-\ep^{-3})^2 (1-\ep^{-2})} \deg \gamma_1-\frac{1}{5}\sum_{\ep \in \mu_5,\ep \neq 1} \frac{\ep^{d-7}}{(1-\ep^{-7})^2 (1-\ep^{-3})}\deg \gamma_2,
\end{eqnarray*}
where $d\deg H|_{C}$ is given by $c_1(\Oh_{X}(d)|_{C})$, and $\gamma_1,\gamma_2$ are the first Chern classes of $\Oh_{C}(3)$, $\Oh_{C}(7)$ respectively. Moreover, we know that the canonical class of $C$ is given by $c_1(\Oh_{C}(10))$ and $\deg H|_C=5 \cdot \frac{2}{125}$. 

Then the remaining parts come from these $3$ singular points, and hence they are given by:
\begin{eqnarray*}
r_3&=&\frac{1}{3}\sum_{\ep\in\mu_3,\ep \neq 1} \frac{\ep^d}{(1-\ep^{-1})^3}+\frac{1}{7}\sum_{\ep\in \mu _7}\frac{\ep^d}{(1-\ep^{-4})(1-\ep^{-5})^2}+\\
&&\frac{1}{25}\sum_{\ep
\in \mu_{25},\ep^5\neq
1}\frac{\ep^d}{(1-\ep^{-3})(1-\ep^{-7})(1-\ep^{-15})}.
\end{eqnarray*}
Using the Program \ref{contribution}, we can calculate the Dedekind sums in the formula. The following are codes in MAGMA program. 
\begin{verbatim}
>h:=2/2625;
>r1:=func<d|(1/6*d^3+341/2*d)*h>;
>s1:=Contribution(5,[2,3]);
>s2:=Contribution(5,[3,3,2]);
>s3:=Contribution(5,[3,2,2]);
>kc:=10; ga1:=3; ga2:=7;
>r2:=func<d|(s1[d mod 5+1]*d-1/2*kc*s1[d mod 5+1]-ga1*s2[(d-3) 
mod 5+1] -ga2*s3[(d-7) mod 5+1])*2/25>;
>c1:=Contribution(3,[1,1,1]);
>c2:=Contribution(7,[4,5,5]);
>c3:=Contribution(25,[3,7,15]);
>r3:=func<d|c1[d mod 3+1]+c2[d mod 7+1]+c3[d mod 25+1]>;
>rr:=[r1(d)+r2(d)+r3(d): d in [2..10]];
>rr;
[ 0, 1, 0, 1, 1, 1, 1, 1, 2 ]
\end{verbatim}
The last output gives the plurigenera for degree $2,\dots, 10$.

\end{exa}

\section{Parcing Hilbert series}\label{Parcing}

\subsection{Statement of the theorem and some examples}
 Let $X$ be a quasismooth variety of dimension $n$. Suppose $X$ has a basket of orbifold curves $\sB_{C}=\{\text{curves of type }\frac{1}{r}(a_1,\dots,a_{n-1})\}$ and a basket of orbifold points $\sB_P=\{\text{points of type }\frac{1}{s}(b_1,\dots,b_n)\}$. 

Using the formula in Proposition \ref{RRwithcurvesingularities}, we can write the Hilbert series associated to $(X,H)$ into the following form as we did for the isolated case in \cite{BRZ}.
\begin{prop}\label{roughdec}
Let $X$ be a quasismooth projective orbifold with polarization $\Oh(1)$. Let $\sB_p$ and $\sB_{C}$ be the orbifold loci given above. Then the Hilbert series $P(t)=\sum_{d\geq 0} h^0(\Oh(d))t^d$ can be written as 
\[
P(t)=\frac{A(t)}{(1-t)^{n+1}}+\sum_{Q\in \sB_P}P_{\mathrm{per},Q}(t)+ \sum_{C\in \sB_{C}} P_{\mathrm{per},C}(t),
\]
where $A(t)$ is a polynomial of degree $k_{X}+n+1$ if $k_{X}\geq 0$; otherwise $A(t)$ is of degree $n$. The term $P_{\mathrm{per}}(t)$ for a point $Q$ of type $\frac{1}{s}(b_1,\dots,b_n)$ is given by
\[
P_{\mathrm{per},Q}(t)=\frac{\sum_{i=1}^{s-1}\frac{1}{s}\sum_{\ep \in \mu_s,\ep^{b_i} \neq 1}\frac{\ep^i}{(1-\ep^{-b_1})\cdots (1-\ep^{-b_{n}})}t^i}{1-t^s},
\]
and the term $P_{\mathrm{per},\sC}(t)$ for a curve $\sC$ of type $\frac{1}{r}(a_1,\dots,a_{n-1})$ is given by
\begin{eqnarray*}
P_{\mathrm{per},\sC}(t)&=&\frac{\sum_{i=1}^{r}i\sigma_it^i}{1-t^r} \deg H|_{\sC}+ \frac{(\sum_{i=1}^{r}\sigma_it^i)t^r}{(1-t^r)^2} r \deg H|_{\sC}-\\
&&\frac{\sum_{i=0}^{r-1}\sigma_i t^i}{1-t^r}\frac{1}{2}\deg K_{X}|_{\sC}-\sum_{j=1}^{n-1}\frac{\sum_{i=0}^{r-1}\delta_{i,j}t^i}{1-t^r}\frac{1}{2}\deg \gamma_j, 
\end{eqnarray*}
where $\sigma_i=\sigma_i(\frac{1}{r}(a_1,\dots,a_{n-1}))$ is given by $\frac{1}{r}\sum_{\ep \in \mu_r,\ep \neq 1} \frac{\ep^i}{(1-\ep^{-a_1})\cdots(1-\ep^{-a_{n-1}})}$ and $\delta_{i,j}=\frac{1}{r}\sum_{\ep\in \mu_r,\ep \neq 1}\frac{\ep^{i}(1+\ep^{-a_j})}{(1-\ep^{-a_j})^2\prod_{i\neq j}(1-t^{a_i})}$. The $\gamma_i$ are given as before.
\end{prop}

\begin{pf}
To see this, just note that the first term in the Riemann--Roch formula is a polynomial in $d$ of degree $n$ and the contributions from points are periodic. Also note that 
\begin{eqnarray*}
&&\frac{a_1t+2a_2t^2+\cdots+r a_{r}t^{r}}{1-t^r}+\frac{(a_1t+a_2t^2+\cdots+a_rt^r) rt^r}{(1-t^r)^2}\\
&&=a_1t+2a_2t^2+\cdots+(r-1)a_{r-1}t^{r-1}+r a_rt^r+(r+1)a_1t^{r+1}+\\
&&(r+2)a_2t^{r+2}+\cdots.
\end{eqnarray*}
For more details, see the proof in Section 3.2 in \cite{BRZ}. $\square$
\end{pf}

\smallskip

The above parsing roughly gives us how each orbifold locus appears in the Hilbert series, but we want a parsing with each of the parts corresponding to orbi\-fold loci characterized in a closed form, analogue to Theorem 1.3 in \cite{BRZ}. The following theorem parses the Hilbert series in such a way.
\begin{theorem}\label{curveparsing}
Let $X$ be a quasismooth projective variety of dimension $n$ with a polarization $\Oh(1)$. Suppose $(X,\Oh(1))$ is projectively Gorenstein, and $X$ has a basket of orbifold curves $\sB_{C}= \{\text{curve } C$ $\text{of type } \frac{1}{r}(a_1,\dots,a_{n-1})\}$ and a basket of orbifold points $\sB_Q=\{\text{point }Q \\\text{ of type } \frac{1}{s}(b_1,\dots, b_n)\}$. Then the Hilbert series associated to $(X,\Oh(1))$ can be uniquely parsed into the form
\[
P(t)=P_I(t)+\sum_{Q \in \sB_Q}P_{\mathrm{orb},\, Q}(t)+\sum_{C\in \sB_{C}} P_{\mathrm{orb},\,C}(t),
\]
where 
\begin{enumerate}
  \item  the initial term $P_I(t)$ is of the form $\frac{I(t)}{(1-t)^{n+1}}$, where $I(t)$ is a polynomial of degree $c=k_{X}+n+1$ and palindromic. $P_I(t)$ has the same coefficients as $P(t)$ as power series up to and including degree $\rd{\frac{c}{2}}$. 
  \item  the orbifold term $P_{\mathrm{orb},Q}(t)$ for a point $Q$ of type $\frac{1}{s}(b_1,\dots,b_n)$ is given by $\frac{Q(t)}{(1-t)^n h(1-t^s)}$, where $h=\mathrm{gcd}\,((1-t^{b_1})\cdots(1-t^{b_n}), \frac{1-t^s}{1-t})$ and $Q(t)$ is the inverse of $\prod \frac{1-t^{b_i}}{1-t}$ modulo $\frac{1-t^s}{(1-t)h}$ supported in $[\rd{\frac{c}{2}}+1+\deg h , \rd{\frac{c}{2}}+s-1]$. For each $Q$, the numerator $Q(t)$ has integral coefficients and $P_{\mathrm{orb}, \,Q}(t)$ is Gorenstein symmetric of degree $k_{X}$.  
  \item the orbifold term $P_{\mathrm{orb},\,C}(t)$ for a curve $C$ of type $\frac{1}{r}(a_1,\dots, a_{n-1})$ can be given in two parts, that is, 
\begin{equation}\label{orbifoldcurvepart}
g_{C}(t)\frac{S_{C,1}(t)}{(1-t)^{n-1}(1-t^r)^2}+\frac{S_{C,2}(t)}{(1-t)^{n}(1-t^r)},
\end{equation}
where 
\begin{itemize}
\item $S_{C,1}(t)$ is given by the inverse of $\prod_{i=1}^{n-1} \frac{1-t^{a_i}}{1-t}$ mod $\frac{1-t^r}{1-t}$, supported in the integral $[\rd{\frac{c+r}{2}}+1,\rd{\frac{c+r}{2}}+r-1]$. Then $S_{C,1}(t)$ has integral coefficients and $\frac{S_{C,1}(t)}{(1-t)^{n-1}(1-t^r)^2}$ is Gorenstein symmetric of degree $k_{X}$.
\item $g_{C}(t)$ is a Laurent polynomial with integral coefficients, which is supported in $[-\rd{\frac{r}{2}}+1,-\rd{\frac{r}{2}}+r-1]$, and $g_{C}(t)$ is palindromic centered at degree $0$. Moreover, $g_{C}(t)$ is determined by the degree of the curve and the dissident points it passes through, as described in Section \ref{withdissidentpoints}. In particular, when there are no dissident points on $C$, $g(t)=r\deg H|_{C}$ is an integer.
\item $\frac{S_{C, 2}(t)}{(1-t)^{n}(1-t^r)}$ has integral coefficients and is Gorenstein symmetric of degree $k_{X}$.
\end{itemize}
\end{enumerate}
\end{theorem}

\bigskip
\begin{rem}
The point of this theorem is to state explicitly how each term is constructed from orbipoints, orbicurves, their normal bundle and the global canonical weight. However, to give a complete description of $S_{C,2}(t)$ in terms of the normal bundle of the curve is still work in progress. 
\end{rem}

We will prove this theorem in the later sections step by step. Now we want to give some examples to verify (or clarify) the statements in the theorem. 
\begin{exa}\label{nondissidentcurve} 
Let $X_{12}$ be a general degree $12$ hypersurface inside $\underline{\PP}^4(1,2,2,3,4)$ with polarization $\Oh(1)$. Then $k_{X}=0$ and $c=0+3+1=4$.
Note that it has an orbicurve $C_{12}\subset \underline{\PP}(2,2,4)$ of degree $3/2$ of type $\frac{1}{2}(1,1)$. The Hilbert series associated to $(X, \Oh_{X}(1))$ can be parsed into
\begin{eqnarray*}
P(t)=\frac{1-t^{12}}{(1-t)(1-t^2)^2(1-t^3)(1-t^4)}= P_I(t)+P_{C}(t)
\end{eqnarray*}
where 
\begin{itemize}
\item  $P_I(t)=\frac{1-3t+5t^2-3t^3+t^4}{(1-t)^4}$ is the initial term. Written as power series, $P_I(t)=1+t+3t^2+7t^3+\cdots$ while $P(t)=1+t+3t^2+4t^3+\cdots$. 
\item $P_{C}(t)=3\frac{-t^3}{(1-t)^2(1-t^2)^2}$. Here we do not have the second part in (\ref{orbifoldcurvepart})of the orbifold curve term  (see a general statement in Proposition \ref{halfcurveterm}). The coefficient $3$ is given by $2\deg H|_{C}$ because there is no dissident points on the curve.  
\end{itemize}
\end{exa}

\begin{exa}\label{anexampleofthedissidentpoint}
Take a general hypersurface $X$ of degree $36$ inside ${\PP}^5(1,4,5,6,9,10)$. We can analyze the orbifold loci on $X$. It has two types of orbifold points, namely the point $P_1=(0,\dots,0,1)$ of type $\frac{1}{10}(1,4,5,9)$ and $2$ points $P_2,P_3$ on the coordinate axis $x_0=x_1=x_2=x_5=0$ of type $\frac{1}{3}(1,1,1,2)$. The $P_1$ is a dissident point, and it lives on the curve $C=C_{36} \subset {\PP}(4,6,10)$ of type $\frac{1}{2}(1,1,1)$ as well as the curve $L=\underline{\PP}(5,10)$ of type $\frac{1}{5}(1,4,4)$. Given the polarization $\Oh(1)$ on $X$, the associated Hilbert series $P(t)$ can be parsed into
\begin{eqnarray*}
P(t)&=&\frac{1-t^{36}}{(1-t)(1-t^4)(1-t^5)(1-t^6)(1-t^9)(1-t^{10})}  \\
&=&P_I(t)+P_{\mathrm{orb}, \,P_1}(t)+ P_{\mathrm{orb}, \, P_2}(t)+P_{\mathrm{orb}, \, P_3}(t)+ P_{\mathrm{orb}, C}(t)+P_{\mathrm{orb},L}(t),
\end{eqnarray*}
where 
\begin{itemize}
\item the intial term $P_I(t)=\frac{1-4t+6t^2-4t^3+6t^4-4t^5+t^6}{(1-t)^5}$.
\item the orbifold point terms are given by $P_{\mathrm{orb}, \,P_1}(t)=\frac{-t^9+t^{10}-t^{11}}{(1-t)^2(1-t^2)(1-t^5)(1-t^{10})}$, and $P_{\mathrm{orb}, \, P_2}(t)=P_{\mathrm{orb}, \, P_3}(t)=\frac{-t^4}{(1-t)^4(1-t^3)}$.
\item the orbifold curve term $P_{\mathrm{orb},C}(t)=0 \frac{-t^4}{(1-t)^3(1-t^2)^2}$ and the orbifold curve term $P_{\mathrm{orb},L}(t)=(t+1/t)\frac{t^7}{(1-t)^3(1-t^5)^2}+\frac{-2t^4-3t^5-2t^6}{(1-t)^4(1-t^5)}$.
\end{itemize}
Note that the degree of the curve $C=C_{36} \subset \underline{\PP}(4,6,10)$ is $3/10$, but gives no contribution in this parsing. This is because the dissident point $P_1$ ``bites off'' its contribution $3/5\frac{-t^4}{(1-t)^3(1-t^2)^2}$. Similarly, for the curve $L=\underline{\PP}(5,10)$, which is of degree $1/10$, the dissident point $P_1$ ``bites off'' $(-t+1/2-1/t)\frac{t^7}{(1-t)^3(1-t^5)^2}$ from this curve contribution, and $g_L(t)$ is given by $5 \deg H|_{C}-(-t+1/2-1/t)=(t+	1/t)$. We will explain what ``bite off'' means in Section \ref{withdissidentpoints}. 
\end{exa}

\subsection{Contributions from dissident points}\label{diss}
Now we start a proof of Theorem \ref{curveparsing}. We consider the formula in Proposition \ref{roughdec} piece by piece and try to adjust each of them to be of the form described in our theorem. Note that the parts corresponding to isolated orbifold points can be treated in the same way as in Theorem $1.3$ in \cite{BRZ}, so we only need to consider the remaining parts, namely the parts corresponding to orbifold curves and dissident points. This section deals with the contribution from dissident points. 

For an orbifold point of type $\frac{1}{s}(b_1,\dots,b_n)$, \textit{dissident} means that there exists some $b_i$ such that $(s,b_i) \neq 1$. Furthermore, if we assume that the orbifolds we consider here only have orbifold loci of dimension $\leq 1$, then there do not exist $i,\,j$ such that $\mathrm{gcd}\,(s,b_i,b_j) \neq 1$. In this case, for each of the $w_i=\mathrm{gcd}\,(s,b_i)\neq 1$, there is a curve of type $\frac{1}{w_i}(\overline{b_1},\dots,\widehat{b_i},\dots,\overline{b_n})$ passing through this point, where $\widehat{b_i}$ means that $b_i$ is omitted, and $\overline{b_j}$ gives the smallest nonnegative residue of $b_j$ mod $w_i$. 

Recall that the periodic term from a dissident point $Q$ of type $\frac{1}{s}(b_1,\dots,b_n)$ in the Hilbert series is given by
\[
P_{\mathrm{per},Q}(t)=\frac{\sum_{i=0}^{s-1}\frac{1}{s}\sum_{\ep \in \mu_s,\ep^{b_i}\neq 1} \frac{\ep^{i}}{(1-\ep^{-b_1})\cdots(1-\ep^{-b_n})}t^i 
}{1-t^s}. 
\]
By Proposition \ref{relation}, the numerator of $P_{\mathrm{per},Q}(t)$, denoted by $N_{\mathrm{per},Q}(t)$,  satisfies
\begin{eqnarray}\label{inversemodrelation}
N_{\mathrm{per}}(t)\prod_{i=0}^n (1-t^{b_i})=1 \text{\, mod \, } \frac{1-t^s}{(1-t) h},
\end{eqnarray}
where $h=\text{gcd }(\prod_{i=1}^n(1-t^{b_i}), \frac{1-t^s}{1-t})$. As in \cite{BRZ}, we want to move some other parts in the Hilbert series to $P_{\mathrm{per},Q}(t)$ so that we obtain $P_{\mathrm{orb},Q}(t)$ with integral coefficients and satisfying the Gorenstein symmetric property. 
\begin{lem}\label{facts}
\begin{enumerate}
	\item Let $f(t)\in \QQ[t]$ be a palindromic polynomial supported in $[\gamma+1, \gamma+l]$ with $0\leq l\leq r-1$. Then given $m \in \ZZ$, there is a unique polynomial $g(t)=f(t) \text{ mod } \frac{1-t^r}{1-t}$ supported in $[\gamma+mr+1, \gamma+(m+1)r-1]$, and obviously $g(t)$ is also palindromic.
	\item If $f(t) \in \QQ[t]$ is palindromic, supported in $[\gamma+1, \gamma+l-1]$, then there exists a palindromic polynomial $g(t)=f(t) \text{ mod }\frac{1-t^r}{1-t}$ with support in $[\gamma+\rd{\frac{l}{2}}+2,\gamma+ \rd{ \frac{l}{2}}+r] $ when $l$ is odd, and with support in  $[\gamma+\frac{l}{2}+2, \gamma+\frac{l}{2}+r-1] $ when $l$ is even.
\end{enumerate}
\end{lem}
\begin{pf}
For the first part, it is easy to see that we only need to shift the degree of each term up or down by $|mr|$. For the second part, just note that subtracting $a_1t^{\gamma+1}\frac{1-t^r}{1-t}$ from $f(t)$ will cancel out two terms, namely $a_1t^{\gamma+1}+a_{1}t^{\gamma+l}$, and do the similar process to the resulting polynomial. We will finally obtain a palindromic polynomial with support as stated. $\square$
\end{pf}

\begin{prop}\label{dissident}
Let $w_i=\mathrm{gcd }\,(s,b_i)$. There exists a unique $Q(t)$ supported in $[\rd{\frac{c}{2}}+1+\deg h, \rd{\frac{c}{2}}+s-1]$ given by the equation
\[
\frac{Q(t)}{\prod_{i=1}^n(1-t^{w_i})(1-t^s)}=P_{\mathrm{per},Q}(t)+\frac{A(t)}{(1-t)^{n+1}}+\sum_{1\leq i \leq n,\, w_i\neq 1}\frac{B_i(t)}{(1-t^{w_i})^2}
\]
where $A(t), \,B_i(t)$ are some Laurent polynomials, and $Q(t)$ can be determined by
\[
Q(t)\prod_{i=1}^n \frac{1-t^{b_i}}{1-t^{w_i}}= 1 \,\mathrm{mod}\, \frac{1-t^s}{(1-t)h},
\]
that is, $Q(t)$ is the inverse of $\prod_{i=1}^{n}\frac{1-t^{b_i}}{1-t^{w_i}}\,\mathrm{mod}\,\frac{1-t^s}{(1-t)h}$. Furthermore, $Q(t)$ has integral coefficients, and $\frac{Q(t)}{(1-t)^nh(1-t^s)}$, denoted by $P_{\mathrm{orb}, Q}(t)$, is Gorenstein symmetric of degree $k_{X}$.  
\end{prop}

\begin{pf}
Note that the equality can be rewritten as 
\[
Q(t)=N_{\mathrm{per},Q}(t)(1-t)^nh+A(t)\frac{(1-t^s)h}{1-t}+\sum_{1\leq i \leq n,\, w_i\neq 1 }B_i(t)(1-t)^n\frac{(1-t^s)h}{(1-t^{w_i})^2},
\]
and in our case $\prod_{i=1}^n\frac{1-t^{w_i}}{1-t}=h$. Therefore, one can write the above equality as
\[
Q(t)=N_{\mathrm{per},Q}(t)\prod_{i=1}^{n}(1-t^{w_i})+\frac{1-t^s}{(1-t)h}(A(t)h^2+\sum_{1\leq i \leq n,\, w_i \neq 1} B_i(t) \frac{(1-t)^{n+1}h^2}{(1-t^{w_i})^2}).
\]
By the above equality and (\ref{inversemodrelation}), we deduce that $Q(t)$ is the inverse of $\prod\frac{1-t^{b_i}}{1-t^{w_i}}$ mod $\frac{1-t^s}{(1-t)h}$.

Moreover, suppose $w_{i_1},\dots,w_{i_k}$ are all the $w_i$ that are not equal to $1$. Then $h=\text{gcd }(\prod_{i=1}^n(1-t^{b_i}), \frac{1-t^s}{1-t})=\prod_{j=1}^k \frac{1-t^{w_{i_j}}}{1-t}$ and gcd $(h^2,\frac{(1-t)^{n+1}h^2}{(1-t^{w_{i_1}})^2},\dots, \frac{(1-t)^{n+1}h^2}{(1-t^{w_{i_k}})^2})=1$. Then by the same idea as in Theorem $2.2$ in \cite{BRZ}, there is a unique $Q(t)$ supported in  $[\rd{\frac{c}{2}}+1+\deg h, \rd{\frac{c}{2}}+s-1]$.

To see that $Q(t)$ has integral coefficients, note that the inverse of $\prod_{i=1}^n \frac{1-t^{b_i}}{1-t^{w_i}}$ can be given by $\prod_{i=1}^n \frac{1-t^{\alpha_ib_i}}{1-t^{b_i}} $ mod $\frac{1-t^s}{(1-t)h}$, where $\alpha_i$ is the smallest positive integer such that $\alpha_ib_i=w_i$ mod $s$. Since $Q(t)$ with length $\leq s-\deg h-1$ can be obtained by moving $\prod_{i=1}^n \frac{1-t^{\alpha_ib_i}}{1-t^{b_i}}$ modulo $\frac{1-t^s}{(1-t)h}$, we conclude that $Q(t)$ has integral coefficients.  

To prove the Gorenstein symmetry of $P_{\mathrm{orb},Q}(t)$, we reduce the support of the polynomial $\prod_{i=1}^n \frac{1-t^{\alpha_ib_i}}{1-t^{b_i}}$ modulo $\frac{1-t^s}{1-t}$ and then modulo $\frac{1-t^s}{(1-t)h}$. Since $\prod_{i=1}^n \frac{1-t^{\alpha_ib_i}}{1-t^{b_i}}$ and $\frac{1-t^r}{(1-t)h}$ as polynomials are both palindromic, we can prove that for the chosen support of $Q(t)$, the orbifold term $P_{\mathrm{orb},Q}(t)$ is Gorenstein symmetric of degree $k_{X}$. Here we show one of the cases, and the rest are similar. 
 
Note that $\prod_{i=1}^n \frac{1-t^{\alpha_ib_i}}{1-t^{b_i}}$ is a polynomial of degree $\sum_{i=1}^n (\alpha_i-1)b_i$ and $\frac{1-t^s}{1-t}$ is a polynomial of degree $s-1$. Suppose $\sum_{i=1}^n (\alpha_i-1)b_i +1$ and $s$ are both even. Then by trimming $\prod_{i=1}^n \frac{1-t^{\alpha_ib_i}}{1-t^{b_i}}$ modulo $\frac{1-t^s}{1-t}$ from both ends, we obtain a palindromic polynomial of length $s-2$ supported in
\begin{equation}\label{(1)}
[\frac{\sum_{i=1}^n(\alpha_i-1)b_i-1}{2}-\frac{s-2}{2}+1, \frac{\sum_{i=1}^n(\alpha_i-1)b_i-1}{2}+\frac{s-2}{2}],
\end{equation}
and by moving a bit forward (see Lemma \ref{facts}, $2$) we can also get another palindromic polynomial supported in 
\begin{equation}\label{(2)}
[\frac{\sum_{i=1}^n(\alpha_i-1)b_i+1}{2}+1, \frac{\sum_{i=1}^n(\alpha_i+1)b_i-1}{2}+s-2].
\end{equation}
Then we trim them further modulo $\frac{1-t^s}{(1-t)h}$. If $\deg h$ is even, then we obtain from (\ref{(1)}) a palindromic polynomial supported in 
\[
[\frac{\sum_{i=1}^n(\alpha_i-1)b_i-1}{2}-\frac{s-2-\deg h}{2}+1, \frac{\sum_{i=1}^n(\alpha_i-1)b_i-1}{2}+\frac{s-2-\deg h}{2}],
\]
and we obtain from (\ref{(2)}) a palindromic polynomial supported in 
\[
[\frac{\sum_{i=1}^n(\alpha_i-1)b_i+1}{2}+\frac{\deg h}{2}+1, \frac{\sum_{i=1}^n(\alpha_i+1)b_i-1}{2}+\frac{\deg h}{2}+(s-\deg h -1)-1].
\]
Notice that 
\begin{eqnarray*}
&&\sum_{i=1}^n(\alpha_i-1)b_i-1-(s-2-\deg h)\\
&=&\sum_{i=1}^n \alpha_ib_i-\sum_{i=1}^nb_i-1-s+\deg h +2\\
&=&\sum_{i=1}^n (w_i-1)+n-\sum_{i=1}^nb_i-1-s+\deg h +2 \,\mathrm{mod}\, s \\
&=&2\deg h +n+k_{\sX}+1 \,\mathrm{mod}\,s,
\end{eqnarray*}
since we know that $\sum_{i=1}^n b_i+k_{\sX}=0 \,\mathrm{mod}\,r$ and $\deg h=\sum_{i=1}^n (w_i-1)$. Therefore we can finally use Lemma \ref{facts}, $1$ to move the support to 
\[
[\frac{c}{2}+\deg h +1, \frac{c}{2}+s-2].
\] 
If $\deg h$ is odd, we just need to replace $\frac{\deg h}{2}$ by $\rd{\frac{\deg h}{2}}$ and replace$\frac{c}{2}$ by $\rd{\frac{c}{2}}$, and the rest of arguments are similar. 
  
Thus we obtain in the end a palindromic polynomial with integral coefficients supported in $[\rd{\frac{c}{2}}+\deg h+1, \rd{\frac{c}{2}}+s-1]$ which is the inverse of $\prod_{i=1}^n \frac{1-t^{b_i}}{1-t^{w_i}}$ mod $\frac{1-t^s}{(1-t)h}$.  $\square$
\end{pf}

\bigskip

\begin{rem}
We should remark here that we made a choice of the form for the dissident point contribution in our Hilbert series parsing. This choice gives us integral coefficients for the numerator of $P_{\mathrm{orb}}(t)$, but it also gives us the denominator of $P_{\mathrm{orb}}$ in the form $(1-t^{w_1})\cdots(1-t^{w_n})(1-t^s)$ for a dissident point of type $\frac{1}{s}(b_1,\dots, b_n)$, where $w_i=\text{gcd } (b_i, s)$. Using this choice, to obtain $P_{\mathrm{orb}}(t)$ we have to move some parts of the terms to $P_{\mathrm{per}}(t)$ from curves that pass through this point as well as some growing part (see Section \ref{withdissidentpoints}).
\end{rem}
\begin{rem}\label{moreprecise}
We have a more precise description of the support of the palindromic polynomial $Q(t)$, that is, when the coindex $c=k_{\sX}+n+1$ is even, the support of $Q(t)$ is in $[\rd{\frac{c}{2}}+\deg h +1, \rd{\frac{c}{2}}+s-2]$; when the coindex $c=k_{\sX}+n+1$ is odd, the support of $Q(t)$ is in $[\rd{\frac{c}{2}}+1, \rd{\frac{c}{2}}+s-1]$. 
\end{rem}
\begin{rem}
Notice that $\prod \frac{1-t^{b_i}}{1-t}/h$ and $\frac{1-t^s}{1-t}/h$ have no common factors. Hence, we can calculate $Q(t)$ using the $\mathrm{XGCD}$ in the MAGMA program, i.e., the inverse of $\prod \frac{1-t^{b_i}}{1-t}/h$ mod $\frac{1-t^s}{1-t}/h$ is given by $\alpha(t)$ in the following equality:
\[
\alpha(t)\prod \frac{1-t^{b_i}}{1-t}/h+\beta(t) \frac{1-t^s}{1-t}/h =1,
\]
and one can shift the support of $\alpha(t)$ to get $Q(t)$.  The following program is analogue to Program in \cite{BRZ}, but it applies to a wider range of types of orbifold points (including the isolated case), that is, it applies to dissident points on curves or dissident points on a higher dimensional orbifold locus. The following program is obtained with help of M. Reid.

\begin{pro}\label{Qorb}
\begin{verbatim}
function Qorb(r,LL,k)
L := [Integers() | i : i in LL];
if (k + &+L) mod r ne 0
   then error "Error: Canonical weight not compatible";
end if;
n := #L; Pi := Denom(L);
A := (1-t^r) div (1-t); B := Pi div (1-t)^n;
H := GCD(A, B); M := &* [GCD(A, 1-t^i) : i in L];
shift := Floor(Degree(M*H)/2);
l := Floor((k+n+1)/2+shift+1);
de := Maximum(0,Ceiling(-l/r));
m := l + de*r;
G, al_throwaway, be := XGCD(A div H, t^m*B div M);
return t^m*be/(M*(1-t)^n*(1-t^r)*t^(de*r));
end function;
\end{verbatim}
\end{pro}
\end{rem}

\subsection{Contributions from curves}\label{curv}
This section deals with the parts that correspond to orbicurves in our parsing. Recall from Proposition \ref{roughdec} that for an orbicurve of type $\frac{1}{r}(a_1,\dots, a_{n-1})$, the original shape of its contribution to the Hilbert series is given by the following:
\begin{eqnarray}\label{curvecontributioninroughdec}
P_{C}(t)&=&\frac{\sum_{i=1}^{r}i\sigma_it^i}{1-t^r} \deg H|_{C}+ \frac{(\sum_{i=1}^{r}\sigma_it^i)t^r}{(1-t^r)^2} r \deg H|_{C}-\\
&&\frac{\sum_{i=0}^{r-1}\sigma_i t^i}{1-t^r}\frac{1}{2}k_{X}\deg H|_{C}-\sum_{j=1}^{n-1}\frac{\sum_{i=0}^{r-1}\delta_{i,j}t^i}{1-t^r}\frac{1}{2}\deg \gamma_j,
\end{eqnarray}
where $\sigma_i=\sigma_i(\frac{1}{r}(a_1,\dots,a_{n-1}))$ is given by $\frac{1}{r}\sum_{\ep \in \mu_r,\ep \neq 1} \frac{\ep^i}{(1-\ep^{-a_1})\cdots(1-\ep^{-a_{n-1}})}$ and $\delta_{i,j}=\frac{1}{r}\sum_{\ep\in \mu_r,\ep \neq 1}\frac{\ep^{i}(1+\ep^{-a_j})}{(1-\ep^{-a_j})^2\prod_{i\neq j}(1-\ep^{-a_i})}$.

\bigskip

We want to show that the above expression can be adjusted to the form $\frac{M(t)}{(1-t)^{n-1}(1-t^r)^2}$, which is Gorenstein symmetric of degree $k_{X}$. We first deal with the parts related to the normal bundle, namely, $\frac{\sum_{i=0}^{r-1}\delta_{i,j}t^i}{1-t^r}\frac{1}{2}\deg \gamma_j$, for $1\leq j\leq n-1$.
\begin{lem}\label{normalbundlepart}
There exists a unique $N_j(t)$ supported in $[\rd{\frac{c}{2}}+1,\rd{\frac{c}{2}}+r-1]$ in the following:
\[
\frac{N_j(t)}{(1-t)^n(1-t^r)}=\frac{\sum_{i=0}^{r-1}\delta_{i,j}t^i}{1-t^r}+\frac{A_j(t)}{(1-t)^{n+1}},
\]
for each $1\leq j\leq n-1$. Moreover, $N_j(t)$ satisfies
\[
N_j(t)\frac{1-t^{a_1}}{1-t}\cdots (\frac{1-t^{a_j}}{1-t})^2 \cdots \frac{1-t^{a_{n-1}}}{1-t}=1+t^{a_j} \,\mathrm{mod}\, \frac{1-t^r}{1-t}
\]
for all $j$. Consequently, $N_j(t)$ has integral coefficients and $\frac{N_j(t)}{(1-t)^n(1-t^r)}$ is Gorenstein symmetric of degree $k_{X}$.
\end{lem}
\begin{pf}
Observe that 
\[
(\sum_{i=0}^{r-1}\delta_{i,j}t^i)(1-t^a_1)\cdots (1-t^{a_j})^2\cdots (1-t^{a_{n-1}})=1+t^{a_j} \,\mathrm{mod}\, \frac{1-t^{r}}{1-t}.
\]
Then the rest follows as we did before. $\square$
\end{pf}

Now we are going to study the first three terms in (\ref{curvecontributioninroughdec}). Putting these three terms together we have 
\[
\frac{\sum_{i=0}^{r-1} (i-\frac{k_{X}}{2})\sigma_{i} t^i+ \sum_{i=0}^{r-1}(r-i+\frac{k_{X}}{2})\sigma_i t^{r+i}}{(1-t^r)^2}\deg H_{C}.
\]
By adding some growing term we can write this into the following form
\[
\frac{N(t)}{(1-t)^{n-1}(1-t^r)^2}=\frac{\sum_{i=0}^{r-1} (i-\frac{k_{X}}{2})\sigma_{i} t^i+ \sum_{i=0}^{r-1}(r-i+\frac{k_{X}}{2})\sigma_i t^{r+i}}{(1-t^r)^2}+\frac{V(t)}{(1-t)^{n+1}},
\]
where $N(t)$ is supported in $[\rd{\frac{c}{2}}+1,\rd{\frac{c}{2}}+2r-2]$. Therefore $N(t)$ is given by 
\begin{eqnarray}\label{xxx}
(\sum_{i=0}^{r-1} (i-\frac{k_{X}}{2})\sigma_{i} t^i+ \sum_{i=0}^{r-1}(r-i+\frac{k_{X}}{2})\sigma_it^{r+i})(1-t)^{n-1} 
\end{eqnarray}
moved to the right support modulo $(\frac{1-t^r}{1-t})^2$.
\begin{lem}\label{mainpartofthecurve}
$N(t)$ is a palindromic polynomial and $\frac{N(t)}{(1-t)^{n-1}(1-t^r)^2}$ is Gorenstein symmetric of degree $k_{X}$. 
\end{lem}
\begin{pf}
To prove that $N(t)$ is palindromic, the idea is that we first move the support of the polynomial (\ref{xxx}) to $[\rd{\frac{c}{2}}, \rd{\frac{c}{2}}+2r-1]$ modulo $(1-t^r)^2$ and then move the support to $[\rd{\frac{c}{2}}+1,\rd{\frac{c}{2}}+2r-2]$ modulo  $(\frac{1-t^r}{1-t})^2$. Note that for any integer $b$ we have
\begin{eqnarray*}
&&t^b(\sum_{i=0}^{r-1} (i-\frac{k_{X}}{2})\sigma_{i} t^i+ \sum_{i=0}^{r-1}(r-i+\frac{k_{X}}{2})\sigma_it^{r+i})\\
&=&\sum_{i=0}^{r-1}(-b+i-\frac{k_{X}}{2})\sigma_{-b+i}t^i+\sum_{i=0}^{r-1}(r+b-i+\frac{k_{X}}{2})\sigma_{-b+i}t^{r+i}\,\mathrm{mod}\,(1-t^r)^2.
\end{eqnarray*}
Now using this equality, for any integer $\gamma$ we obtain
\begin{eqnarray*}
&&t^\gamma(1-t)^{n-1}(\sum_{i=0}^{r-1} (i-\frac{k_{X}}{2})\sigma_{i} t^i+ \sum_{i=0}^{r-1}(r-i+\frac{k_{X}}{2})\sigma_it^{r+i})\\
&\equiv& \sum_{j=0}^{n-1}(-1)^j {n-1 \choose j} t^{\gamma+j}(\sum_{i=0}^{r-1} (i-\frac{k_{X}}{2})\sigma_{i} t^i+ \sum_{i=0}^{r-1}(r-i+\frac{k_{X}}{2})\sigma_it^{r+i})\\
&\equiv&\sum_{i=0}^{r-1} \sum_{j=0}^{n-1}(-1)^j{n-1 \choose j} (-(\gamma+j)+i-\frac{k_{X}}{2}) \sigma_{-(\gamma+j)+i}t^i+\\
&&\sum_{i=0}^{r-1}\sum_{j=0}^{n-1}(-1)^j{n-1 \choose j }(r+(\gamma+j)-i+\frac{k_{X}}{2})\sigma_{-(\gamma+j)+i}t^{r+i},
\end{eqnarray*}
where $\equiv$ means equality modulo $(1-t^r)^2$. Here we want to show that if we choose $\gamma=-\rd{\frac{c}{2}}$, the last polynomial above is palindromic, and we denoted it by $L_{\gamma}(t)$. 

Now let $\rho_i$ be the coefficient of degree $i$ in $L_{\gamma}(t)$. We show that when $c$ is even, $\rho_{2r-1}=0$ and $\rho_{i}=\rho_{2r-2-i}$; when $c$ is odd, $\rho_{i}=\rho_{2r-1-i}$. Therefore $L_{\gamma}(t)$ is palindromic in the support $[0,\rd{\frac{c}{2}}+2r-1]$. Here we only show it for the case when $c$ is even; the other case is similar. When $c$ is even, we have
\begin{eqnarray*}
\rho_{2r-1}&=&\sum_{j=0}^{n-1}(-1)^j{n-1 \choose j }(r+(-\frac{c}{2}+j)-(r-1)+\frac{k_{X}}{2})\sigma_{-(-\frac{c}{2}+j)+r-1}\\
&=&\sum_{j=0}^{n-1}(-1)^j{n-1 \choose j}j\sigma_{-\frac{c}{2}-j-1}+(\frac{1-n}{2})\sum_{j=0}^{n-1}(-1)^j{n-1 \choose j} \sigma_{\frac{c}{2}-j-1}\\
&=&\frac{1}{r}\sum_{\ep \in \mu_r,\ep \neq 1} \frac{\ep^{\frac{c}{2}-1}(\sum_{j=0}^{n-1}(-1)^j{n-1 \choose j}j \ep^{-j})}{(1-\ep^{-a_1})\cdots (1-\ep^{-a_{n-1}})}+\\
&&\frac{1-n}{2}\frac{1}{r}\sum_{\ep\in \mu_r,\ep \neq 1}\frac{\ep^{\frac{c}{2}-1}(\sum_{j=0}^{n-1}(-1)^j{n-1 \choose j} \ep^{-j})}{(1-\ep^{-a_1})\cdots (1-\ep^{-a_{n-1}})}\\
&=&(1-n)\frac{1}{r}\sum_{\ep \in \mu_r,\ep \neq 1} \frac{\ep^{\frac{c}{2}-2}(1-\ep^{-1})^{n-2}}{(1-\ep^{-a_1})\cdots (1-\ep^{-a_{n-1}})}+\\
&&\frac{1-n}{2}\frac{1}{r}\sum_{\ep\in \mu_r,\ep \neq 1}\frac{\ep^{\frac{c}{2}-1}(1-\ep^{-1})^{n-1}}{(1-\ep^{-a_1})\cdots (1-\ep^{-a_{n-1}})}\\
&=& \frac{1-n}{2}\frac{1}{r}\sum_{\ep \in \mu_r,\ep\neq 1} \frac{\ep^{\frac{c}{2}-1}(1+\ep^{-1})(1-\ep^{-1})^{n-2}}{(1-\ep^{-a_1})\cdots (1-\ep^{-a_{n-1}})}.
\end{eqnarray*}
By the above expression, we can see that $\rho_{2r-1}=0$ by the following fact:
\begin{eqnarray*}
\rho_{2r-1}&=& \frac{1-n}{2}\frac{1}{r}\sum_{\ep \in \mu_r,\ep\neq 1} \frac{\ep^{-\frac{c}{2}+1}(1+\ep)(1-\ep)^{n-2}}{(1-\ep^{a_1})\cdots (1-\ep^{a_{n-1}})}\\
&=&\frac{1-n}{2}\frac{1}{r}\sum_{\ep \in \mu_r,\ep\neq 1} \frac{(-1)^{n-2}\ep^{-\frac{c}{2}+1+n-1}(1+\ep^{-1})(1-\ep^{-1})^{n-2}}{ (-1)^{n-1}\ep^{a_1+\cdots+a_{n-1}}(1-\ep^{-a_1})\cdots (1-\ep^{-a_{n-1}})}\\
&=&- \rho_{2r-1},
\end{eqnarray*}
where the last equality is due to the fact that $\sum_{i=1}^{n-1}a_i=-k_{X} $ mod $r$. 

To prove $\rho_i=\rho_{2r-2-i}$, we first simplify the expression of $\rho_i$ and $\rho_{2r-2-i}$. First, we have
\begin{eqnarray*}
\rho_i&=&\sum_{j=0}^{n-1}(-1)^j{n-1 \choose j} (-(-\frac{c}{2}+j)+i-\frac{k_{X}}{2}) \sigma_{-(-\frac{c}{2}+j)+i}\\
&=& \sum_{j=0}^{n-1}(-1)^j{n-1 \choose j}(-j+\frac{n+1}{2}+i) \sigma_{-j+\frac{c}{2}+i}\\
&=&\sum_{j=0}^{n-1}(-1)^j{n-1 \choose j}(-j)\sigma_{-j+\frac{c}{2}+i}+\sum_{j=0}^{n-1}(-1)^j{n-1\choose j}\sigma_{-j+\frac{c}{2}+i}(\frac{n+1}{2}+i)\\
&=&\frac{1}{r}\sum_{\ep\in \mu_r,\ep \neq 1}\frac{(n-1)(1-\ep^{-1})^{n-2}\ep^{\frac{c}{2}+i-1}}{(1-\ep^{-a_1})\cdots(1-\ep^{-a_{n-1}})}+\\
&&\frac{1}{r}\sum_{\ep\in \mu_r,\ep \neq 1}\frac{(1-\ep^{-1})^{n-1}\ep^{\frac{c}{2}+i}}{(1-\ep^{-a_1})\cdots(1-\ep^{-a_{n-1}})}(\frac{n+1}{2}+i),
\end{eqnarray*}
where the last equality uses the binomial expansion $(1-t)^{n-1}=\sum_{j=0}^{n-1}(-1)^j{n-1\choose j}t^j$ and its derivative $(n-1)t(1-t)^{n-2}=\sum_{j=0}^{n-1}(-1)^j j{n-1 \choose j}t^{j}$.


Recall that we have $\sigma_{\frac{c}{2}-l}=(-1)^{n-1} \sigma_{\frac{c}{2}+l-n-1}$. Therefore we can simplify $\rho_{2r-2-i}$ as follows:
\begin{eqnarray*}
\rho_{2r-2-i}&=&\sum_{j=0}^{n-1}(-1)^j{n-1 \choose j }(r+(-\frac{c}{2}+j)-(r-2-i)+\frac{k_{X}}{2})\sigma_{-(-\frac{c}{2}+j)+r-2-i}\\
&=&\sum_{j=0}^{n-1}(-1)^j{n-1 \choose j }(j+\frac{-n+3}{2}+i)\sigma_{\frac{c}{2}-i-j-2}\\
&=&\sum_{j=0}^{n-1}(-1)^j{n-1 \choose j }(-(n-1-j)+\frac{n+1}{2}+i)(-1)^{n-1}\sigma_{\frac{c}{2}+i+j+2-n-1}\\
&=&(-1)^n\sum_{j=0}^{n-1}(-1)^j{n-1 \choose j }(n-1-j)\sigma_{\frac{c}{2}+i-((n-1)-j)}+\\
&&(-1)^{n-1}\sum_{j=0}^{n-1}(-1)^j{n-1 \choose j }\sigma_{\frac{c}{2}+i-((n-1)-j)}(\frac{n+1}{2}+i)\\
&=&\frac{1}{r}\sum_{\ep\in \mu_r,\ep \neq 1}\frac{(n-1)(1-\ep^{-1})^{n-2}\ep^{\frac{c}{2}+i-1}}{(1-\ep^{-a_1})\cdots(1-\ep^{-a_{n-1}})}+\\
&&\frac{1}{r}\sum_{\ep\in \mu_r,\ep \neq 1}\frac{(1-\ep^{-1})^{n-1}\ep^{\frac{c}{2}+i}}{(1-\ep^{-a_1})\cdots(1-\ep^{-a_{n-1}})}(\frac{n+1}{2}+i),
\end{eqnarray*}
where the last equality uses the fact that $\sum_{j=0}^{n-1}(-1)^j(n-1-j)\ep^{-(n-1-j)}=(n-1)\ep^{-1}(\ep^{-1}-1)^{n-2}$. 

Therefore, $\rho_i=\rho_{2r-2-i}$ for all $0\leq i \leq r-1$. Since $N(t)$ supported in $[\rd{\frac{c}{2}}+1, \rd{\frac{c}{2}}+2r-2]$ is given by $L_{-\rd{\frac{c}{2}}}$ modulo the palindromic polynomial $(\frac{1-t^r}{1-t})^2$, we obtain that $N(t)$ is also palindromic. $\square$ 
\end{pf}

\smallskip

Combining Lemma \ref{normalbundlepart} and Lemma \ref{mainpartofthecurve}, we can get the following conclusion, which gives a global view of the curve contribution in our Hilbert series parsing.
\begin{prop}
The total contribution from a curve $C$ of type $\frac{1}{r}(a_1,\dots,a_{n-1})$ is given by
\[
\frac{M(t)}{(1-t)^{n-1}(1-t^r)}=\frac{N(t)}{(1-t)^{n-1}(1-t^r)^2}\deg H|_{C} +\sum_{j=1}^{n-1}\frac{N_j(t)}{(1-t)^{n}(1-t^r)} \frac{\deg \gamma_i}{2},
\]
and we this denote by $P_{C}(t)$. Moreover, $P_{C}(t)$ is Gorenstein symmetric of degree $k_{X}$.
\end{prop}

\bigskip

Even though we can prove that the curve contribution $P_{C}(t)$ has the Gorenstein symmetry property, we cannot characterize it as a whole using some ice cream function (or InverseMod function) as we did for point contributions.  However, we can give a characterization as an ice cream function for the ``order $2$'' part, more precisely, we can put the part 
\[
P_{\mathrm{per},\,C}(t)=\frac{t^r\sum_{i=1}^{r}\frac{1}{r}\sum_{\ep \in \mu_r,\ep \neq 1}\frac{\ep^i}{(1-\ep^{-a_1})\cdots(1-\ep^{-a_{n-1}})}t^i}{(1-t^r)^2}
\]
into an ice cream function by the following lemma.
\begin{lem}\label{curvepartcontribution}
There exists a unique $S_1(t)$ supported in $[\rd{\frac{c+r-1}{2}}+1,\rd{\frac{c+r-1}{2}}+r-1]$, satisfying
\[
\frac{S_1(t)}{(1-t)^{n-1}(1-t^r)^2}=P_{\mathrm{per},\,C}(t)+\frac{B(t)}{1-t^r}+\frac{A(t)}{(1-t)^{n+1}}.
\]
Consequently, $S_1(t)$ can be determined by the inverse of $\prod \frac{1-t^{a_i}}{1-t}\,\mathrm{mod}\, \frac{1-t^r}{1-t}$ with the chosen support. Moreover, $S_1(t)$ has integral coefficients and $\frac{S_1(t)}{(1-t)^{n-1}(1-t^r)^2}$ is Gorenstein symmetric of degree $k_{X}$.
\end{lem}
\begin{pf}
See proof for isolated case in \cite{BRZ}. $\square$
\end{pf}

Taking into consideration of the coefficient of $P_{\mathrm{per},C}(t)$, we know that 
\[
P_{C}(t)- r\deg H|_{C} \frac{S_1(t)}{(1-t)^{n-1}(1-t^r)^2}
\]
is of the form $\frac{S_2(t)}{(1-t)^n(1-t^r)}$, which is also Gorenstein symmetric of degree $k_{X}$. That is,
\[
P_{C}(t)=r\deg H|_{C} \frac{S_1(t)}{(1-t)^{n-1}(1-t^r)^2}+\frac{S_2(t)}{(1-t)^n(1-t^r)},
\]
and we denote the first part by $P_{C,1}(t)$ and the second part by $P_{C,2}(t)$.

For curves without dissident points this gives a nice form, since $r\deg H|_{C}$ is an integer (see proof of Proposition \ref{integral}). When there are no dissident points on the curve $C$, then $P_{C}(t)$ gives us the $P_{\mathrm{orb},C}(t)$ in our theorem \ref{curveparsing}. However, when there are dissident points on the curve, the number $r\deg H|_{C}$ is possibly fractional. We will see in the next section how orbifold terms we chose for the dissident points affect the number $r\deg H|_{C}$.

\subsection{Orbicurves with dissident points}\label{withdissidentpoints}
Recall that in Proposition \ref{dissident}, choosing $P_{\mathrm{orb},Q}(t)$ in our parsing, we need to move some parts from the curve terms in the Hilbert series to $P_{\mathrm{per},\,Q}(t)$. We will see that after subtracting all the parts which the dissident points ``bite off'', the remaining curve contributions have integral coefficients. Here we can only measure how much the dissident point ``bites off'' from the first part of the curve contribution $P_{C,1}(t)$. We cannot control precisely how the dissident points affect the second part $P_{C,2}(t)$, but we can prove what each dissident point "bites off" from the second part is Gorenstein symmetric of degree $k_{X}$. 
 
\begin{prop}\label{biteoff}
Let $Q$ be a dissident point of type $\frac{1}{s}(b_1,\dots,b_n)$. Let $w_i=(s,b_i)$ and $P_{\mathrm{orb},Q}(t)=\frac{Q(t)}{(1-t^{w_1})\cdots (1-t^{w_n})(1-t^s)}$ be the term given in Proposition \ref{dissident}. Then when $w_i\neq 1$, there is a curve $C_i$ of type $\frac{1}{w_i}(\overline{b_1},\dots, \widehat{b_i},\dots,\overline{b_n})$ passing through this point. Then the point $Q$ bites off the following contribution from $P_{C_i,1}(t)$:
\[
\mathrm{bit}_{Q,w_i}(t) \frac{S_{1,w_i}(t)}{(1-t)^{n-1}(1-t^{w_i})^2},
\]
where $S_{1,w_i}(t)$ is given as in Lemma \ref{curvepartcontribution}. The coefficient $\mathrm{bit}_{Q,w_i}(t)$ is a Laurent polynomial supported in $[-\rd{\frac{w_i}{2}}+1, \rd{\frac{w_i}{2}}-1]$ and is Gorenstein symmetric of degree $0$ $($in the sense that $\mathrm{bit}_{Q,w_i}(t)=(t)^0\mathrm{bit}_{Q,w_i}(1/t)$$)$, determined uniquely by 
\[
\mathrm{bit}_{Q,w_i}(t)=\frac{w_i}{s} Q(t)  \prod_{j\neq i}\frac{1-t^{b_j}}{1-t^{w_j}}\,\,\mathrm{mod}\,\, \frac{1-t^{w_i}}{1-t},
\] 
Moreover, $\mathrm{bit}_{Q,w_i}(t)$ has integral coefficients except for the constant term. 
\end{prop}

\begin{pf}
Recall from Proposition \ref{dissident} that the orbifold term $P_{\mathrm{orb}}(t)$ for a dissident point $Q$ of type $\frac{1}{s}(b_1,\dots, b_{n})$ is given by
\[
\frac{Q(t)}{\prod_{i=1}^n(1-t^{w_i})(1-t^s)}=\frac{N_{\mathrm{per},Q}(t)}{1-t^s}+\frac{A(t)}{(1-t)^{n+1}}+\sum_{1\leq i \leq n, w_i\neq 1}\frac{B_i(t)}{(1-t^{w_i})^2}.
\]
We can rewrite this in the form
\begin{eqnarray*}
&&\frac{Q(t)}{\prod_{i=1}^n(1-t^{w_i})(1-t^s)}=\frac{N_{\mathrm{per},Q}(t)}{1-t^s}+\frac{A^\prime(t)}{(1-t)^{n+1}}\\
&&+ \sum_{1\leq i \leq n, w_i \neq 1} (\mathrm{bit}_{Q,w_i}(t)\frac{S_{1,w_i}(t)}{(1-t)^{n-1}(1-t^{w_i})^2}
+\frac{D_i(t)}{(1-t)^n(1-t^{w_i})}),
\end{eqnarray*}
which gives
\begin{eqnarray*}
Q(t)&=&N_{\mathrm{per},Q}(t)\prod_{i=1}^{n}(1-t^{w_i})+ \frac{1-t^s}{(1-t)h}(A'(t)h^2+\\
&& \sum_{1\leq i \leq n, w_i\neq 1}(\mathrm{bit}_{Q,w_i}S_{1,w_i}(t)+D_i(t)\frac{1-t^{w_i}}{1-t})(\prod_{j\neq i}\frac{1-t^{w_j}}{1-t})^2),
\end{eqnarray*}
where $h=$ gcd $(\prod_{i=1}^n(1-t^{b_i}), \frac{1-t^s}{1-t})=\prod_{i=1}^n\frac{1-t^{w_i}}{1-t}$ and each $S_{1,w_i}(t)$ is the inverse of $\prod_{j\neq i}\frac{1-t^{b_j}}{1-t}$ mod $\frac{1-t^{w_i}}{1-t}$. Note that $h^2$ and $(\prod_{j\neq i}\frac{1-t^{w_j}}{1-t})^2$, for $i=1,\dots,n$, are coprime, which ensures that we can move $Q(t)$ to the right support; $S_{1,w_i}(t)$ and $\frac{1-t^{w_i}}{1-t}$ are coprime, which enables us to choose $\mathrm{bit}_{Q,w_i}(t)$ modulo $\frac{1-t^{w_i}}{1-t}$.  We claim that we can choose $\mathrm{bit}_{Q,w_i}(t)$ to be Gorenstein symmetric of degree $0$. In fact, by the above equality, we know that $\mathrm{bit}_{Q,w_i}$ satisfies 
\begin{eqnarray*}
\mathrm{bit}_{Q,w_i}(t)&\equiv&Q(t)\frac{(1-t)h}{1-t^s}(\prod_{j\neq i}\frac{1-t}{1-t^{w_j}})^2\prod_{j\neq i}\frac{1-t^{b_j}}{1-t}\\
&\equiv&Q(t)\frac{1}{1+t^{w_i}+\cdots+t^{\frac{s}{w_i}-1}} \prod_{j\neq i}\frac{1-t^{b_j}}{1-t^{w_j}}\\
&\equiv&\frac{w_i}{s}Q(t) \prod_{j\neq i}\frac{1-t^{b_j}}{1-t^{w_j}}\,\, \mathrm{mod}\,\, \frac{1-t^{w_i}}{1-t}.
\end{eqnarray*}
Since $Q(t)$ and $\prod_{j\neq i}\frac{1-t^{b_j}}{1-t^{w_j}}$ are symmetric, we deduce that $\mathrm{bit}_{Q,w_i}(t)$ can be reduced to be Gorenstein symmetric of degree $0$ modulo $ \frac{1-t^{w_i}}{1-t}$ (This can be done as in Section \ref{anotherproof}). Moreover, we know that the constant part of $\mathrm{bit}_{Q,w_i}$ is given by $\frac{\alpha_iw_i}{s}$ plus an integer, where $\alpha_i$ is the smallest positive integer such that $\alpha_i b_i= w_i$ mod $s$, and apart from the constant term, $\mathrm{bit}_{Q,w_i}$ has integral coefficients. In fact, recall that $Q(t)=\prod_{i=1}^n \frac{1-t^{\alpha_i b_i}}{1-t^{b_i}}+\beta(t)\frac{1-t^s}{(1-t)h}$. We plug this into the above equality and get 
\begin{eqnarray*}
&&\mathrm{bit}_{Q,w_i}(t)\equiv \frac{w_i}{s}(\prod_{i=1}^n \frac{1-t^{\alpha_i b_i}}{1-t^{b_i}}+\beta(t)\frac{1-t^s}{(1-t)h}) \prod_{j\neq i}\frac{1-t^{b_j}}{1-t^{w_j}}\\
&&\equiv(\frac{w_i}{s}\frac{1-t^{\alpha_ib_i}}{1-t^{b_i}}\prod_{j\neq i}\frac{1-t^{\alpha_ib_j}}{1-t^{w_j}}+\frac{w_i}{s}\beta(t)\frac{1-t^s}{1-t^{w_i}}\prod_{j\neq i}\frac{1-t^{\beta_jw_j}}{1-t^{w_j}}\prod_{j\neq i}\frac{1-t^{b_j}}{1-t^{w_j}})\\
&&\equiv (\frac{\alpha_iw_i}{s}+\beta(t)\prod_{j\neq i}\frac{1-t^{\beta_jw_j}}{1-t^{w_j}}\prod_{j\neq i}\frac{1-t^{b_j}}{1-t^{w_j}}) \,\,\mathrm{mod}\,\,\frac{1-t^{w_i}}{1-t},
\end{eqnarray*}
where $\beta_j$ satisfies $\beta_jw_j=1$ mod $w_i$. Such $\beta_j$ exist because $w_j$ and $w_i$ are coprime for $j\neq i$. Note that the second part of the last equality, $\beta(t)\prod_{j\neq i}\frac{1-t^{\beta_jw_j}}{1-t^{w_j}}\prod_{j\neq i}\frac{1-t^{b_j}}{1-t^{w_j}}$, is a polynomial with integral coefficients. Since $\mathrm{bit}_{Q,w_i}(t)$ is uniquely determined with chosen support, then the constant term of $\mathrm{bit}_{Q, w_i}(t)$ is given by $\frac{\alpha_iw_i}{s}$ plus some integer, and apart from the constant term, it has only integral coefficients.  $\square$
\end{pf}

\bigskip
Here we give one example to explain the last proposition.
\begin{exa}
Given a point $Q$ of type $\frac{1}{10}(1,4,5,9)$, we have $w_1=w_4=1$, $w_2=2$ and $w_3=5$. Then $Q$ lies on both a curve of type $\frac{1}{2}(1,1,1)$ and a curve of type $\frac{1}{5}(1,4,4)$. By the last proposition, it bites off $\mathrm{bit}_{Q,w_2}$ from the curve of type $\frac{1}{2}(1,1,1)$, which is given by $3/5\frac{S_{1,w_2}(t)}{(1-t)^3(1-t^2)^2}$. In fact, $P_{\mathrm{orb},Q}(t)$ can be calculated using Program \ref{Dorb}, which gives us
\[
P_{\mathrm{orb},Q}(t)=\frac{-t^9+t^{10}-t^{11}}{(1-t)^2(1-t^2)(1-t^5)(1-t^{10})}.
\]
By the above proposition, we know that
\[
\mathrm{\mathrm{bit}}_{Q,w_2}(t)=\frac{2}{10}(-t^9+t^{10}-t^{11})\frac{1-t^9}{1-t} \,\mathrm{mod}\, \frac{1-t^2}{1-t}=3/5.
\]
Similarly, for $\mathrm{bit}_{Q,w_3}(t)$ we have
\[
\mathrm{bit}_{Q,w_3}(t)=\frac{5}{10}(-t^9+t^{10}-t^{11}) \frac{1-t^4}{1-t^2}\frac{1-t^9}{1-t } \,\mathrm{mod}\,\frac{1-t^5}{1-t}= -t+1/2-1/t.
\]
\end{exa}

\begin{rem}
Note that the parts, which a dissident point bites off from each of the curves it lies on,  are determined by its orbifold type and do not depend on the ambient orbifold it lives in.
\end{rem}
\smallskip
Now we know how each dissident point affects the curves it lies on. Given a curve $C$ of type $\frac{1}{r}(a_1,\dots, a_{n-1})$ with a set $\sT$ of dissident points on it, we have the following:
\begin{prop}\label{integral}
$r\deg H|_{C}-\sum_{Q \in \sT} \mathrm{bit}_{Q,r}(t)$ has integral coefficients and is Gorenstein symmetric of degree $0$.
\end{prop} 
\begin{pf}
The only thing we need to prove is that its constant term is an integer. Since we only consider the case when there are only orbifold loci of dimension $\leq 1$, for each dissident point $Q\in \sT$ of type $\frac{1}{s_Q}(b_{Q,1},\dots, b_{Q,n})$, there exists exactly one $b_{Q,i}$ such that gcd $(b_{Q,i},s_Q)=r$ and gcd $(b_{Q,j},r)=1$ for all $j\neq i$. For convenience, we denote this $b_{Q,i}$ by $b_{Q}$.  Recall that the constant term each dissident point bites off from the curve contribution is given by $\frac{\alpha_{Q}r}{s_Q}$ plus some integer, where $\alpha_{Q} b_{Q}=r $ mod $ s_Q$, which is also equivalent to $\alpha_{Q} \frac{b_{Q}}{r}=1$ mod $\frac{s_Q}{r}$.

Since this only concerns the curve $C$, we can restrict the problem to $C$. Suppose $C$ is defined by $I$ in $\underline{\PP}(c_1,\dots,c_l)$, where the $c_i$ are divisible by $r$.  Consider the curve $C^\prime$ defined by the same ideal $I$ in $\underline{\PP}(\frac{c_1}{r},\dots, \frac{c_l}{r})$ with shifted weights for each of the variables. Then the degree $\deg H^\prime|_{C^\prime}$ of the curve $C^\prime$ is given by $r\deg H|_{C}$. The dissident point $Q$ restricted to the curve $C^\prime$ is an orbifold point of type $\frac{1}{s_Q/r}(\frac{b_Q}{r})$. Recall in Section \ref{nonwellformed} that the Euler characteristic of $\Oh_{C^\prime}(1)$ is given by $\chi(\Oh_{C^\prime})+\deg H'|_{C^\prime}-\sum_Q \frac{\alpha_Q}{s_Q/r}$ if the curve has orbifold points of type $\frac{1}{s_Q/r}(b_Q/r)$. Thus, we see that $r\deg H|_{C}-\sum_Q \frac{\alpha_Q r}{s_Q}=\deg H^\prime|_{C^\prime}-\sum_{Q} \frac{\alpha_Q}{s_Q/r}=\chi(\Oh_{C^\prime}(1))-\chi(\Oh_{C^\prime})$ is an integer. We are done.  $\square$
\end{pf}

\begin{exa}
Now we can return to Example \ref{anexampleofthedissidentpoint} to work out the coefficients for the $P_{C,1}(t)$ for each of the curves.
\end{exa}

\bigskip

So far, we have only considered how dissident points on the curve affects the first part $P_{C,1}(t)$ of the curve term. Now we want to see how the second term $P_{C,2}(t)$ is affected by the dissident points. Note that even though we cannot control precisely the parts that the dissident points bite off from the second piece $P_{C,2}(t)$, we can assert the following:
\begin{prop}\label{secondpartofthecurve}
Let $Q$ be an orbifold point of type $\frac{1}{s}(b_1,\dots, b_n)$ on $X$. Suppose $w_i=\,\mathrm{ gcd }\,(s,b_i)\neq 1$ for $1\leq i\leq l$ (possibly after reordering the $b_i$), and let $C_i$ be the orbifold curve of type $\frac{1}{w_i}(b_1,\dots, \widehat{b_i},\dots, b_n)$ that passes through $Q$.  Then with the $P_{\mathrm{orb},Q}(t)$ given in Proposition \ref{inversemodrelation}, $P_{\mathrm{orb},Q}(t)$ ``bites off'' from the second part $P_{C_i,2}$ a rational function that is Gorenstein symmetric of degree $k_{X}$.
\end{prop}
\begin{pf}
Note that the numerator $N_{\mathrm{per},\,Q}(t)$ of the periodic term $P_{\mathrm{per},\,Q} (t)$ is divisible by $h(t)=\prod_{i=1}^n\frac{1-t^{w_i}}{1-t}$, where $w_i=\text{ gcd }(b_i,s)$ for all $i$. Then there exists a unique $n(t)$ supported in $[\rd{\frac{c}{2}}+1+\rd{\frac{\deg h}{2}}, \rd{\frac{c}{2}}+r-1-\rd{\frac{\deg h}{2}}]$ in the following equality: 
\[
\frac{n(t)}{(1-t)^n m(t)}=\frac{N_{\mathrm{per},\,Q}(t)}{h(t)m(t)}+\frac{A(t)}{(1-t)^{n+1}},
\]
where $m(t)=\frac{1-t^s}{h(t)}$ and $A(t)$ is some Laurent polynomial. Equivalently, we have
\[
n(t)=\frac{N_{\mathrm{per},\,Q}(t)}{h(t)}(1-t)^n+A(t)\frac{m(t)}{1-t}.
\]
Hence $n(t)$ is the inverse of $\prod_{i=1}^n\frac{1-t^{b_i}}{1-t}$ mod $\frac{1-t^s}{(1-t)h(t)}$ by Proposition \ref{relation}. One can prove that $\frac{n(t)}{(1-t)^n m(t)}$ is Gorenstein symmetric of degree $k_{X}$ as before. Therefore,
\[
P_{\mathrm{orb},Q}(t)-\frac{n(t)}{(1-t)^n m(t)}-\sum_{i=1}^l\mathrm{bit}_{Q,C_i}(t) \frac{S_{1,w_i}(t)}{(1-t)^{n-1}(1-t^{w_i})^2}
\]
is Gorenstein symmetric of degree $k_{X}$, which is of the form \\$\frac{S(t)}{(1-t)^{n+1-l}(1-t^{w_1})\dots(1-t^{w_{l}})}$. This is the sum of what $P_{\mathrm{orb},Q}(t)$ bites off from the second part $P_{C_i,2}(t)$ of each curve $C_i$, that is 
\[
\frac{S(t)}{(1-t)^{n+1-l}(1-t^{w_1})\dots(1-t^{w_{l}})}=\frac{s_1(t)}{(1-t)^n(1-t^{w_1})}+\cdots+\frac{s_l(t)}{(1-t)^n(1-t^{w_l})},
\] 
where $\frac{s_i(t)}{(1-t)^n(1-t^{w_i})}$ represents the bite from the second part of the curve $C_i$, and $s_i(t)$ is supported in $[\rd{\frac{c}{2}}+1, \rd{\frac{c}{2}}+w_i-1]$. Now we need to prove that the Gorenstein symmetry of the sum implies the Gorenstein symmetry of $\frac{s_i(t)}{(1-t)^n(1-t^{w_1})}$ for all $0 \leq i \leq l$. In fact, the above equality can be rewritten as 
\[
S(t)=\sum_{i=1}^l s_i(t) \prod_{j=1,j\neq i}^l \frac{1-t^{w_j}}{1-t}.
\]
Now if we take the last equality modulo $\frac{1-t^{w_i}}{1-t}$, then 
\[
S(t)= s_i(t) \prod_{j=1,j\neq i}^l \frac{1-t^{w_j}}{1-t}\,\mathrm{mod }\, \frac{1-t^{w_i}}{1-t},
\]
which implies that 
\[
s_i(t)=S(t)\prod_{j=1,j\neq i}^l \frac{1-t^{w_j v_{i_j}}}{1-t^{w_j}}\,\mathrm{mod }\, \frac{1-t^{w_i}}{1-t},
\]
where $w_jv_{i_j}=1 $ mod $w_i$. In this way we can prove as before that $s_i(t)$ is Gorenstein symmetric with the support $[\rd{\frac{c}{2}}+1, \rd{\frac{c}{2}}+w_i-1]$. $\square$
\end{pf}

Combining Propositions \ref{integral} and \ref{secondpartofthecurve}, we know that after subtracting what each of the dissident points bites off from the curve, the remaining contribution from the curve $C$ in the Hilbert series is given in the following form:
\begin{equation}\label{orbicurve}
(r\deg H|_{C}-\sum_{Q\in \sT}  \mathrm{bit} _{Q,r}(t) )\frac{S_1(t)}{(1-t)^{n-1}(1-t^r)^2}+ \frac{S_2(t)}{(1-t)^n(1-t^r)},
\end{equation}
where each part is Gorenstein symmetric of degree $k_{X}$. We denote the above expression by $P_{\mathrm{orb},C}(t)$ for a curve with dissident points.

\subsection{A special case}

For an orbifold curve, we have seen that in general its contribution in our Hilbert series parsing consists of two parts as in (\ref{orbicurve}).  The following proposition says that for an orbifold curve of type $\frac{1}{2}(1,\dots,1)$, we only have the first part of the contribution. 
\begin{prop}\label{halfcurveterm}
Let $(X,H)$ be a projectively Gorenstein pair. Suppose there is an orbifold curve of singularity type $\frac{1}{2}(1,\dots,1)$, and that there are  dissident points of type $\sT=\{Q \text{ of type }\frac{1}{2s}(2b_{Q,1},\dots,b_{Q,n})\}$ living on $C$ (by assumption we have $\mathrm{gcd}\,(b_{Q,i},2s)=1$ for all $i$). Then the orbifold term for this curve $C$ can be given by
\begin{equation}
 P_{\mathrm{orb},C}(t)=\alpha \frac{t^{\rd{\frac{c+1}{2}}+1}}{(1-t)^{m-1}(1-t^2)^2},
\end{equation}
where $\alpha=2\deg H|_{C}-\sum_{Q\in \sT} \mathrm{bit}_{Q}(t)$ and $\mathrm{bit}_{Q}(t)$ are determined as in Proposition \ref{biteoff}. 
\end{prop}
\begin{pf}
Note that since $(X, H)$ is projectively Gorenstein, then by Proposition \ref{singularity} we know that $n-1+k=0$ mod $2$. Therefore the coindex $c=k+n+1$ is always even. The second part from the curve contribution is of the form $\frac{t^{\rd{\frac{c}{2}}+1}}{(1-t)^n(1-t^2)}$. When $c$ is even, this part cannot be Gorenstein symmetric of degree $k_{X}$. Then it has to be zero, and so in this case the curve contribution term in our parsing only consists of the first part. $\square$
\end{pf}

\subsection{Initial term and the end of the proof}
Now to finish the parsing of our Hilbert series, we are left with the initial term. Recall that our orbifold $X$ has orbifold curves $\sB_{C}$ and orbifold points $\sB_{Q}$. We have given an orbifold term for each orbifold locus in our parsing of the Hilbert series, namely, $P_{\text{orb},C}(t)$ and $P_{\text{orb},Q }(t)$. Then the remaining part is 
\[
P(t)-\sum_{C\in \sB_{C}} P_{\mathrm{orb},C}(t)-\sum_{Q\in \sB_{Q}}P_{\mathrm{orb},Q }(t),
\]
which we define to be the initial term $P_{I}(t)$. Since each term in the above expression is Gorenstein symmetric of degree $k_{X}$, then $P_I(t)$ is also Gorenstein symmetric of degree $k_{X}$.

Recall that we required the orbifold term in the Hilbert series for points and curves to have numerators with support starting from $\rd{\frac{c}{2}}+1$, and therefore the initial term needs to take care of the first $\rd{\frac{c}{2}}+1$ terms, namely, $P_0,\dots,P_{\rd{\frac{c}{2}}}$, in the Hilbert series. Since $P_I(t)$ is Gorenstein symmetric of degree $k_{X}$, we can write down $P_I(t)$ as in \cite{BRZ}, and $P_I(t)$ has a numerator with integral coefficients by construction. For later use, here we give a MAGMA program to calculate initial term.  That is, given Gorenstein symmetric degree $k$ and the first $\rd{\frac{c}{2}}+1$ initial terms as a vector $L$, we get initial term with the following program.
\begin{pro} \label{initialprogram}
\begin{verbatim}
function initial(L,k,n)
f:=&+[L[i]*t^(i-1): i in [1..#L]];
pp:=R!(f*(1-t)^(n+1));
c:=k+n+1;
if IsEven(c) eq true then 
return (&+[Coefficient(pp, i )*(t^i+t^(c-i)):i in [0..c div 2-1]]+
Coefficient(pp,c div 2)*t^(Floor(c/2)))/(1-t)^(n+1);
else 
return &+[Coefficient(pp,i)*(t^i+t^(c-i)):i in [0..Floor(c/2)]]
/(1-t)^(n+1);
end if;
end function;
\end{verbatim}
\end{pro}

Now we have our parsing as follows:
\[
P(t)=P_I(t)+\sum_{Q\in \sB_Q} P_{\mathrm{orb},Q}(t)+\sum_{C\in \sB_{C}}(P_{C,1}(t)+P_{C,2}(t)).
\]
There is one more point we need to prove, that is, the integral condition for the second part of the curve contribution, namely, $P_{C,2}(t)$ for each orbifold curve $C$. However, we know that the sum $\sum_{C\in\sB_{C}}P_{C,2}(t)$ has integral coefficients, which is of the form
\[
\frac{S(t)}{(1-t)^{n+1}\prod_{C\in \sB_{C}}\frac{1-t^{r_C}}{1-t}}.
\] 
Recall that $P_{C,2}(t)$ is of the form $\frac{S_{C,2}(t)}{(1-t)^n(1-t^{r_{C}})}$. Then
\[
\sum_{C\in\sB_{C}}P_{C,2}(t)=\sum_{C\in \sB_{C}}\frac{S_{C,2}(t)}{(1-t)^n(1-t^{r_{C}})}.
\]
Therefore, $S_{C,2}(t)$ is given by $S(t)(\prod_{C^\prime \neq C}\frac{1-t^{r_{C^\prime}}}{1-t})^{-1}$ mod $\frac{1-t^{r_C}}{1-t}$ (see proof of Proposition \ref{secondpartofthecurve}), which proves that $S_{C,2}(t)$ has integral coefficients as usual. 

This finishes the proof of our Theorem \ref{curveparsing}.

\section{Examples and applications}\label{applications}

In this section, we give some examples of our Hilbert series parsing formula. Then we apply this to construct orbifolds with certain invariants and orbifold loci. First, let us see some examples of our parsing formula with pure orbicurves (that is, orbicurves without dissident points). 
\begin{exa}
Let $X_{10}$ be a degree $10$ hypersurface in $\PP(1,1,1,2,2,2)$ and $\Oh(1)$ be the polarization. This is a canonical $4$-fold with an orbicurve of type $\frac{1}{2}(1,1,1)$. We know $k_{X}=1$ and $c=1+4+1=6$. Also we can calculate the degree of the curve
\[
\deg H|_{C}=\frac{10\cdot 2}{2\cdot 2 \cdot 2}=\frac{5}{2}.
\]
Thus the parsing of Hilbert series is given by 
\[
P(t)=P_I(t)+5 P_{C}(t),
\]
where $P_I$ can be calculated using Program \ref{initialprogram}, which gives
\[
P_I(t)=\mathrm{initial}\,([1,3,9,19],1,4)=\frac{1-2t+4t^2-6t^3+4t^4-2t^5+t^6}{(1-t)^5}.
\]
and $P_{C}$ can be calculated using Program \ref{Qorb}, and it gives
\[
P_{C}(t)=\mathrm{Qorb}\,(2,[1,1,1],3)/(1-t^2)=\frac{t^4}{(1-t)^3(1-t^2)^2}.
\]
\end{exa}

\begin{exa}
Consider the following two $4$-folds:
\begin{itemize}
	\item let $(X_1,\Oh(1))$ be a general hypersurface of degree $16$ in $\PP(1,1,1,3,3,8)$. Then it has an orbicurve $C=\PP(3,3)$ of type $\frac{1}{3}(1,1,2)$. 
	\item let $(X_2, \Oh(1))$ be a general hypersurface of degree $13$ in $\PP(1,1,1,3,3,5)$. Then it has an orbicurve $C^{\prime} =\PP(3,3)$ of type $\frac{1}{3}(1,1,2)$ and an orbipoint of type $\frac{1}{5}(1,1,1,3)$. 
\end{itemize}
Note that these $4$-folds both have canonical weight $-1$ and coindex $c=-1+4+1=4$. They all have the same plurigenera $1,3,6$ in degree $0,1,2$ respectively. Therefore, they have the same initial term, which can be calculated by Program \ref{initialprogram}. This gives
\[
P_{I}(t)=\mathrm{initial}\,([1,3,6],-1,4)=\frac{1-2t+t^2-2t^3+t^4}{(1-t)^5}.
\]
They also both have an orbicurve of type $\frac{1}{3}(1,1,2)$ of the same degree $\frac{1}{3}$, for which we can calculate the first part of the curve contribution by Program \ref{Qorb}, that is, 
\[
P_{C,1}=3\deg H|_{C}\,\mathrm{Qorb}(3,[1,1,2],-1+3)/(1-t^3)=\frac{-t^4}{(1-t)^3(1-t^3)^2}.
\]
where the second part of the curve parsing can be calculated by its Gorenstein property and an extra information of the third plurigenus. Now for $(X_1,\Oh(1))$ we write out our parsing 
\begin{eqnarray*}
P_1(t)&=&P_I(t)+P_{C,1}(t)+P_{C,2}(t)\\
&=&\frac{1-2t+t^2-2t^3+t^4}{(1-t)^5}+\frac{-t^4}{(1-t)^3(1-t^3)^2}+\frac{4t^3}{(1-t)^4(1-t^3)}.
\end{eqnarray*}
For $(X_2,\Oh(1))$, our parsing is
\begin{eqnarray*}
P_1(t)&=& P_I(t)+P_{\mathrm{orb},Q}(t)+P_{C^\prime,1}(t)+P_{C^\prime,2}(t)\\
&=& \frac{1-2t+t^2-2t^3+t^4}{(1-t)^5}+\frac{t^3+t^5}{(1-t)^4(1-t^5)}+\\
&&\frac{-t^4}{(1-t)^3(1-t^3)^2}+\frac{3t^3}{(1-t)^4(1-t^3)}.
\end{eqnarray*}
where $P_{\mathrm{orb},Q}(t)$ is calculated by $\mathrm{Qorb},(5,[1,1,1,3],-1)$ and the second part of the curve contribution is calculated as above.

As one may notice that even though the orbifold types of the two orbicurves $C$ and $C^{\prime}$ are the same, the second parts of the curve contributions are different. This is because the second part of the curve contribution is related to the normal bundle of the curve. 
\end{exa}

Now we have seen some examples of our Hilbert series parsing formula. We want to construct orbifolds with this parsing as in Section \ref{applications}. Here we have a simple example.
\begin{exa}
Suppose we want to construct an orbifold of dimension $3$ with trivial canonical sheaf with the following data:
\begin{itemize}
	\item the first three plurigenera: $P_0=1$, $P_1=1$, $P_2=2$;
	\item an orbicurve $C_1$ of type $\frac{1}{2}(1,1)$ and an orbicurve $C_2$ of type $\frac{1}{3}(1,2)$;
	\item a dissident point $Q_1$ of type $\frac{1}{9}(1,2,6)$ and a dissident point $Q_2$ of type $\frac{1}{6}(1,2,3)$.
\end{itemize}
Suppose such an orbifold exist, then in our Hilbert series parsing we should have 
\[
P_I(t)=\mathrm{initial}\,([1,1,2],0,3)=\frac{1-3t+4t^2-3t^3+t^4}{(1-t)^4}.
\]
We should also have a term related to the curve of type $\frac{1}{2}(1,1)$, that is,
\[
P_{\mathrm{orb},C_1}(t)=\mathrm{Qorb}\,(2,[1,1],2)/(1-t^2)=\frac{-t^3}{(1-t)^2(1-t^2)^2},
\]
and a term related to the curve of type $\frac{1}{3}(1,2)$, which is given by
\begin{eqnarray*}
P_{\mathrm{orb},C_2}(t)&=&P_{C_2,1}(t)+P_{C_2,2}(t)=\mathrm{Qorb}\,(3,[1,2],3)/(1-t^3)+P_{C_2,2}(t)\\
&=&\frac{-t^4}{(1-t)^2(1-t^3)^2}+\frac{S(t)}{(1-t)^3(1-t^3)},
\end{eqnarray*}
where $S(t)$ should be given by $t^3$ multiplied with some integer due to its Gorenstein symmetry property. Moreover, for these two dissident points we should also have orbifold terms 
\begin{eqnarray*} 
P_{\mathrm{orb},Q_1}(t)&=&\mathrm{Qorb}\,(9,[1,2,6],0)=\frac{t^6-t^7+t^8}{(1-t)^2(1-t^3)(1-t^9)};\\
P_{\mathrm{orb},Q_1}(t)&=&\mathrm{Qorb}\,(6,[1,2,3],0)=\frac{t^6}{(1-t)(1-t^2)(1-t^3)(1-t^6)}.
\end{eqnarray*}
To find such an orbifold, we can do the following search:
\begin{verbatim}
pi:=initial([1,1,2],0,3);
q1:=Qorb(2,[1,1],2)/(1-t^2);
q2:=Qorb(3,[1,2],3)/(1-t^3);
q3:=Qorb(9,[1,2,6],0);
q4:=Qorb(6,[1,2,3],0);
for  i,j,k in [0..3] do
p:=pi+i*q1+j*q2 + k*t^3/Denom([1,1,1,3])+q3+q4;
p*Denom([1,2,3,6,9]);[i,j,k];
end for;
\end{verbatim}
Among the outputs (here for simplicity we do not consider the candidates that are codimension $\geq 4$), we have two candidates that possibly gives us such orbifolds, namely, when $i=0$, $j=2$, $k=1$, we have a Hilbert series
\[
P_1(t)=\frac{1- t^9 - 3t^{12}+ 3t^{18} + t^{21} -t^{30}}{(1-t)(1-t^2)(1-t^3)^2(1-t^6)^2(1-t^9)},
\]
and when $i=1$, $j=0$, $k=1$, we have a Hilbert series
\[
P_2(t)=\frac{1- t^{10}- 2t^{12}- t^{13} - t^{15} + t^{16}+ t^{18}+ 2t^{19}+ t^{21}-t^{31}}{(1-t)(1-t^2)(1-t^3)(1-t^4)(1-t^6)^2(1-t^9)}.
\]
Now we analyze these two Hilbert series one by one. In the first case, $P_1(t)$ suggests a codimension $3$ orbifold in $\PP(1,2,3,3,6,6,9)$. Denote the variables of $\PP(1,2,3,3,6,6,9)$ by $x,y,z_1,z_2,t_1,t_2,w$. Then it can be given by $4\times 4$ Pfaffians in the following $5\times 5$ skew symmetric matrix 
\[
\begin{pmatrix}
   &w  &  a_9 & b_9 & c_6 \\
   &   &  t_1 & d_6 & e_3 \\
   &   &      & t_2 & z_1\\
   &   &      &     & z_2
 \end{pmatrix}
\] 
where $a_9$, $b_9$, $c_6$, $d_6$, $e_3$ represent general homogeneous polynomials of degrees $9,9,6,6,3$ respectively. Then the Pfaffians are given by the following equations
\begin{eqnarray*}
pf_1&=&t_1z_2-t_2e_3+z_1d_6,\\
pf_2&=&z_2a_9-b_9z_1+c_6t_2,\\
pf_3&=&wz_2-b_9e_3+d_6c_6,\\
pf_4&=& wz_1-a_9e_3+t_1c_6,\\
pf_5&=&wt_2-a_9d_6+b_9t_1.
\end{eqnarray*}
we can check that the orbifold defined by these equations has the property we required. For example, we see that the point $(0,\dots, 0,1)$ has local parameters $x,y,t_2$, and its orbifold type is given by $\frac{1}{9}(1,2,6)$. Similarly, we can check for other orbifold loci. 

Now in the sencond case, the Hilbert series $P_2(t)$ suggests an orbifold owning these properties can be given by a codimension $3$ orbifold in $\PP(1,2,3,4,6,6,9)$. Denote its coordinates by $(x,y,z,t,w_1,w_2,v)$. Then this orbifold can be defined by Pfaffians in the following matrix
\[
\begin{pmatrix}
   & v &  a_9 & t^2+b_8 & c_6 \\
   &   &    d_7 &  w_2  & y^2+t+e_4 \\
   &   &        &  w_1   &  t\\
   &   &        &       & z 
 \end{pmatrix}
\] 
and we can check that general choices of these homogeneous polynomials will give us an orbifold with the required properties.  
\end{exa}

\end{document}